%% file: main.tex
\documentclass[journal,10pt,letterpaper]{IEEEtran}
\input{preamble}
\begin{document}
\title{Dual-Mode Extremum Seeking for Diffusion Systems with Boundary-Cost Feedback}

\author{Martin Guay$^{1}$ and Denis Dochain$^{2}$ 
	\thanks{$^{1}$Martin Guay is with the Department of Chemical Engineering, Queen's University, 19 Division Street, Kingston, ON, Canada, K7L 3N6. {\tt\small guaym@queensu.ca} 
        $^{2}$ Denis Dochain is with ICTEAM, UCLouvain, B\^atiment Euler, avenue Georges Lema\^itre 4-6, 1348 Louvain-la-Neuve, Belgium {\tt\small denis.dochain@uclouvain.be}}
}

\maketitle
\input{sections/abstract}
\input{sections/introduction}
\input{sections/problem_formulation}
\input{sections/controller}
\input{sections/averaged_stability}
\input{sections/harmonic_lifting}
\input{sections/regularity}
\input{sections/trajectory_comparison}
\input{sections/practical_stability}
\input{sections/frequency_and_scope}
\input{sections/physical_example}
\input{sections/conclusion}
\IEEEtriggeratref{14}
\bibliographystyle{IEEEtran}
\bibliography{references}
\end{document}

%% file: preamble.tex
\usepackage[T1]{fontenc}
\usepackage[utf8]{inputenc}
\usepackage{amsmath,amssymb,amsthm}
\usepackage{microtype,cite}
\usepackage{graphicx,booktabs,array}
\usepackage{xcolor}
\usepackage[hidelinks,bookmarksnumbered=true]{hyperref}
\hypersetup{
 pdftitle={Dual-Mode Extremum Seeking for Diffusion Systems with Boundary-Cost Feedback},
 pdfauthor={Martin Guay},
 pdfsubject={Extremum seeking for boundary-actuated diffusion systems},
 pdfkeywords={extremum seeking, diffusion, boundary control, practical stability}
}
\graphicspath{{figures/}}
\newif\ifmarkup
\markuptrue
\DeclareRobustCommand{\rev}[1]{\ifmarkup{\color{blue}#1}\else{#1}\fi}
\newtheorem{assumption}{Assumption}
\newtheorem{proposition}{Proposition}
\newtheorem{lemma}{Lemma}
\newtheorem{theorem}{Theorem}
\newtheorem{corollary}{Corollary}
\theoremstyle{remark}

\newcommand{\Cb}{\mathcal C_b}
\newcommand{\R}{\mathbb R}
\newcommand{\X}{\mathcal X}
\newcommand{\E}{\mathcal E}
\newcommand{\norm}[1]{\lVert#1\rVert}

\newcommand{\eps}{\varepsilon}
\newcommand{\ee}{\mathrm e}
\newcommand{\dd}{\,\mathrm d}
\newcommand{\Zs}{Z^\star}
\newcommand{\KB}{\mathcal K_B}
\newcommand{\KBref}{\mathcal K_{B,\mathrm{ref}}}
\newcommand{\Romega}{\mathcal R_\omega}
\newcommand{\domega}{\delta_\omega}
\newcommand{\bomega}{\beta_\omega}
\newcommand{\Jk}{\mathcal J_\kappa}

%% file: sections/abstract.tex
\begin{abstract}
{Thermally tuned devices pose a basic control problem. A heater on
one edge of a thin plate must hold a small element at an unknown
optimal temperature, and the only signal available is the device's
own performance.} Boundary measurements of {such} an unknown
performance objective can be used to tune a diffusion process without
reconstructing its spatial state. This paper studies a dual-mode
extremum-seeking controller for {this setting, modeled as} a
semilinear diffusion equation on a square with nonuniform boundary
actuation. The implementation contains two scalar states and exploits
the leading attenuation and \rev{$\pi/4$} phase \rev{lag} of the
diffusion response{, a property shared by every periodically heated
surface}. It requires neither numerical plant coefficients nor the
sign of the input gain. An energy estimate first establishes
exponential stability of a coupled linear diffusion and controller
model. A local bound on the remaining reaction and cost nonlinearities
then yields stability and disturbance estimates for the averaged
system. Complete kernel and boundary-trace arguments justify a change
of coordinates that compares the implemented and averaged trajectories
for all time. The resulting practical stability guarantee concerns the
field, input bias, and sensitivity estimate, and is unaffected by a
constant cost offset under corresponding filter initialization. At
fixed dither amplitude, the residual state bound decreases with the
inverse square root of frequency, while the ultimate cost bound
decreases with its inverse. A thermal optical-filter example
illustrates the predicted trends and the conservatism of the
sufficient conditions.
\end{abstract}

%% file: sections/introduction.tex
\section{Introduction}\label{sec:introduction}
{Many devices perform well only within a narrow window of an
internal variable that nobody measures directly. A thermally tuned
optical filter is a representative case
\cite{Domash2004,HohlfeldZappe2007}. Its transmission peaks at one
temperature of the filter element. This temperature shifts with
fabrication tolerances and ageing. The signal that is convenient
to measure online is the transmitted optical power itself. The
actuator, a resistive heater along one edge of the supporting plate,
is not collocated with the element. The heat it releases must diffuse
across the plate while part of it leaks to the surroundings. The
control task is therefore to drive the heater so that the filter sits
at its transmission peak, through a medium governed by a diffusion
equation, using one scalar performance signal and no thermometer.}

{Extremum-seeking control is the natural tool for a problem of
this kind. It} adjusts a plant input using measured performance when
the optimizing operating condition is unknown. Input dithers reveal local sensitivity, while averaging and
singular perturbation provide classical analysis tools
\cite{KrsticWang2000}. Lie-bracket approximations address a related
class of oscillatory optimization laws \cite{Durr2013}.
Proportional--integral designs combine sensitivity feedback with
input-bias adaptation \cite{Guay2016,GuayDochain2017}. Dual-mode
designs also address the input needed to sustain the operating
condition \cite{GuayAtta2018}, and phasor-based methods exploit
structural information about the plant response \cite{AttaGuay2017}.

Distributed dynamics attenuate and phase-shift the perturbation between
actuation and observation{, and the sensitivity estimate must
account for both}.
Earlier applications include convection--reaction reactors
\cite{Cougnon2006,Hudon2008}. For diffusion actuation, probing-signal
design and boundary compensation have been combined with
infinite-dimensional stability analysis \cite{Feiling2018}. Extensions
address multivariable PDE dynamics \cite{Oliveira2020},
reaction--advection--diffusion equations \cite{Oliveira2021},
spatially integrated outputs \cite{Coutinho2025}, unbiased designs
\cite{Yilmaz2024}, and physics-informed probing \cite{Guo2026}.
Broader developments appear in
\cite{OliveiraKrstic2022,OliveiraKrsticBasar2026}. {These designs
compensate the actuation dynamics online from a known model of the PDE.}

{This paper takes a different route, suggested by the physics of
the thermal problem. When the heater power is modulated sinusoidally,
the temperature at the heated edge responds as a thermal wave. Its amplitude falls with the inverse square root of the frequency and its phase lags the heat flux by $\pi/4$. These are classical properties of diffusion at a periodically heated surface \cite{CarslawJaeger1959}. They hold whatever the conductivity, the heat loss, and the plate size, which enter only through a single amplitude factor carrying the sign of the input gain. Placing the sensor on the actuated edge therefore hands the extremum-seeking loop a response whose phase is known before the controller is switched on, and the demodulation can be fixed once and for all.} These properties permit a controller with two internal states using only the measured
cost. The feedback law requires no numerical PDE coefficients,
state-field measurement, or online PDE compensator. Signal injection
plays a related role in finite-dimensional virtual-measurement
constructions \cite{Combes2016}, but boundary actuation and point
observation require additional trace estimates here.

{A design that uses no model of the plate must be shown to work on
the plate, and} the contribution is a constructive link between
averaged stability and implemented performance. We first prove
exponential stability of the averaged closed loop with its
finite-bandwidth sensitivity filter. A longitudinal reaction
contribution incorporates existing damping, and a local derivative
bound controls the remaining nonlinearities. An exact harmonic lifting,
spatial smoothing, and a boundary-trace argument then give trajectory
closeness for all time. The resulting local practical exponential
regulation has decreasing state and cost bounds as frequency increases,
even at fixed dither amplitude. The comparison and its constants are
developed directly, within the broader context of infinite-dimensional
averaging \cite{HaleVerduynLunel1990} and input-to-state stability
\cite{MironchenkoPrieur2020}.

The objective is regulation toward an attainable local minimum of a
boundary-output cost. The theorem is specific to the stated geometry
and nonlinear margin. The thermal example illustrates controller
performance and spatial regulation, while explicitly distinguishing
its observed frequency trends from the conservative sufficient
certificate.

Sections~\ref{sec:problem}--\ref{sec:controller} specify the problem
and controller. Averaged stability is established in
Section~\ref{sec:stability}. Sections~\ref{sec:harmonic}--\ref{sec:practical}
develop the lifting, regularity, comparison, and practical stability
arguments. Frequency selection and {the thermal filter example}
conclude the analysis. {Readers interested mainly in the design and
its performance can pass from Sections~\ref{sec:problem}--\ref{sec:controller}
directly to Sections~\ref{sec:practical}--\ref{sec:physicalexample}.
Table~\ref{tab:constants} collects the constants that appear in the
final estimates and says where each is defined.}

%% file: sections/problem_formulation.tex
\section{Problem Formulation}\label{sec:problem}
\begin{figure}[t]
\centering
\includegraphics[width=\columnwidth,trim=6 10 7 5,clip]{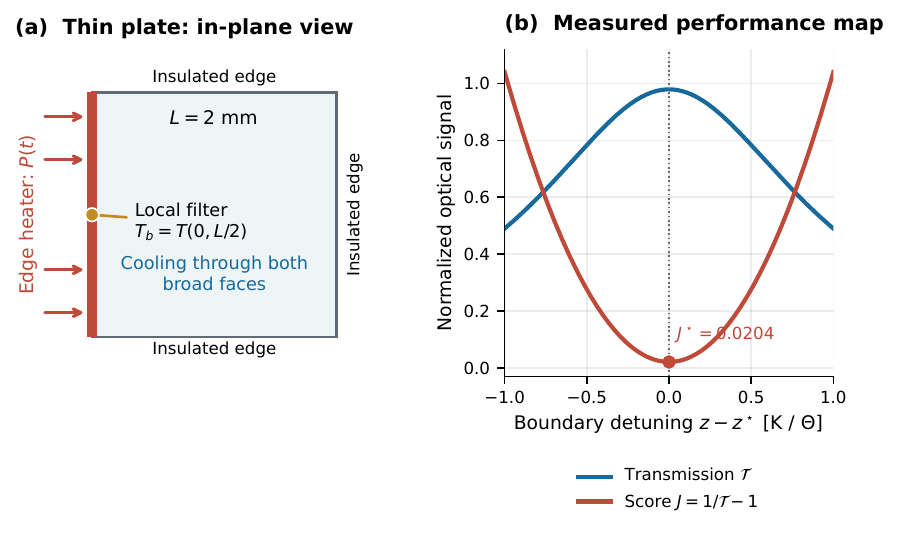}
\caption{{(a) Thin plate heated along one edge, with the optical
element at the midpoint of that edge. The three other edges are
insulated and heat leaks through the two broad faces. (b) Normalized
transmission of the element and the score $J=1/T-1$ against the
detuning of the element temperature from its optimum, in units of the
temperature scale $\Theta$. Only $J$ is measured.}}
\label{fig:geometry}
\end{figure}
{Figure~\ref{fig:geometry}(a) shows the physical situation we have
in mind. A thin square plate carries a small optical element at the
midpoint of one edge. A resistive heater runs along that edge. The
power it delivers may be shaped along the edge but is set by a single
command $u$. Heat spreads through the plate by conduction, leaks to the
surroundings through the two broad faces, and the three remaining
edges are insulated. The transmission of the element depends on its
own temperature and peaks at a value that is not known in advance,
Fig.~\ref{fig:geometry}(b), and the only measurement is the
transmitted optical power, summarized by a score $J$ that is smallest
at the peak. Neither the temperature field nor the temperature of the
element is measured.}

{Scaling lengths by the side of the plate, time by its diffusion
time, and temperature by the width of the transmission peak
 turns the plate into the unit square and the temperature excess $w$ above ambient into the state of a semilinear diffusion equation with boundary actuation. Section~\ref{sec:physicalexample} gives the details.
In this dimensionless form, and allowing for a reaction term and a
fixed heat load that the thermal model does not need but that cost
nothing in the analysis, the plant reads as follows. }
Let $\Omega=(0,1)^2$, $s=(s_1,s_2)$, and consider
\begin{equation}\label{eq:plant}
\begin{aligned}
 w_t&=\kappa\Delta w-\mu_0w+F(s,w)+d(s),\\
 -\kappa w_{s_1}(0,s_2,t)&=bq(s_2)u(t),\\
 w_{s_1}(1,s_2,t)&=0,\\
 w_{s_2}(s_1,0,t)&=w_{s_2}(s_1,1,t)=0,
\end{aligned}
\end{equation}
where $\kappa,\mu_0>0$, $b\ne0$, $u\in\R$, and
$q(s_2)=1+\vartheta\cos(\pi s_2)$ with $\vartheta\in\R$.
The fixed source $d$ may be unknown.
{Here $\kappa$ is the scaled diffusivity, $\mu_0$ accounts for the
heat lost through the faces, $q$ is the distribution of heater power
along the actuated edge, and $b$ is the input gain. The sign of $b$ is unknown such that
the design also covers cooling actuators and sensors of either
polarity.}
Only the noise-free cost
\begin{equation}\label{eq:cost}
 J(t)=\varphi(\Cb w(t)),\qquad \Cb w=w(0,\tfrac12),
\end{equation}
is measured. The field and its boundary value are not separately available.

{Since the heater acts through the boundary and the initial
temperature need not match the boundary condition, the analysis needs
a precise notion of solution.}
Let $\X=C(\overline\Omega;\R)$ be the continuous real-valued fields,
with $\norm f_\infty=\max_{s\in\overline\Omega}|f(s)|$.
The distinct symbol $\Cb:\X\to\R$ denotes point observation,
so $\norm{\Cb}=1$. For $0<\gamma<1$, use
\[
 \norm f_{C^\gamma}=\norm f_\infty+[f]_\gamma,\qquad
 [f]_\gamma=\sup_{s\ne\widetilde s}
 \frac{|f(s)-f(\widetilde s)|}{|s-\widetilde s|^\gamma}.
\]
Let $A=\kappa\Delta-\mu_0$ be the homogeneous Neumann generator,
$S(t)$ its semigroup, and $A_{-1},S_{-1}(t)$ their extrapolations.
Let $B$ be the scalar boundary-command operator in the extrapolation
space. Write $\KB(t)=S_{-1}(t)B$ and
$F(x)=F(\cdot,x(\cdot))$. Solutions have the mild interpretation
\begin{equation}\label{eq:mild}
\begin{aligned}
 w(t)={}&S(t-t_0)w_0
  +\int_{t_0}^t S(t-r)[F(w(r))+d]\,\dd r\\
 &+\int_{t_0}^t\KB(t-r)u(r)\,\dd r.
\end{aligned}
\end{equation}
The kernel estimates below establish local integrability in $\X$.
This formulation permits continuous initial fields without Neumann
compatibility \cite{Nittka2014}.

\begin{assumption}\label{ass:regularity}
The source $d$ is continuous. The functions $F,F_w$ are jointly
continuous near the equilibrium field range and its prescribed
enlargement. There are $H_F,L_F<\infty$ and $0<\gamma<1$ such that
\begin{equation}\label{eq:regularity}
\begin{aligned}
 |F(s,\zeta)-F(\widetilde s,\widetilde\zeta)|
 &\le H_F|s-\widetilde s|^\gamma\\
 &\quad+L_F|\zeta-\widetilde\zeta|
\end{aligned}
\end{equation}
on that region. The cost $\varphi$ is $C^2$ on an open interval
containing the corresponding boundary values.
\end{assumption}

{Assumption~\ref{ass:regularity} asks only for H\"older continuity
in space and Lipschitz continuity in temperature of the reaction, and
for a twice differentiable cost. Both are mild for the thermal
application, where the reaction is absent and the cost is a smooth
function of the element temperature.}

\begin{assumption}\label{ass:optimum}
There is a stationary pair $(w^\star,u^\star)$ of \eqref{eq:plant}
whose observed value $z^\star=\Cb w^\star$ satisfies
$\varphi'(z^\star)=0$ and $h:=\varphi''(z^\star)>0$.
\end{assumption}

{In the thermal system, Assumption~\ref{ass:optimum} says that
some constant heater power holds the element exactly at the temperature
where its transmission peaks, and that the peak is not flat.}
Thus the scalar output minimum is attainable at steady state.
The objective is local regulation of the associated field and bias,
not optimization of an arbitrary steady input-to-cost map.
Neither the equilibrium nor the cost derivatives are known to the
controller. For the nominal estimates, fix $H_0<\infty$ and require
\begin{equation}\label{eq:initialfamily}
 \norm{w_0-w^\star}_{C^\gamma}\le H_0.
\end{equation}
The initial error will also be small in $\X$. This fixed family
controls the integrable startup layer of the measured trace.

{The design relies on one physical observation. A sinusoidal
modulation of the heater power launches a thermal wave into the plate,
and at the frequencies of interest the temperature at the heated edge
responds with an amplitude that decays like $\omega^{-1/2}$ and a
phase that lags the heat flux by $\pi/4$, the classical signature of
diffusion at a periodically heated surface \cite{CarslawJaeger1959}.
Conductivity, heat loss, and plate size do not alter this leading
behavior; they only set the amplitude. In the present notation,}
for the linear diffusion equation with $F=d=0$, the midpoint
transfer function is
\begin{equation}\label{eq:transfer}
 H_b(i\omega)=\frac{b}{\sqrt\kappa\sqrt{\mu_0+i\omega}}
       \coth\sqrt{\frac{\mu_0+i\omega}{\kappa}}.
\end{equation}
{Here the square roots are taken with positive real part, and the
transverse cosine mode of $q$ does not appear because it vanishes at
the midpoint.} Section~\ref{sec:harmonic} proves
\begin{equation}\label{eq:asymptotic}
 H_b(i\omega)=\chi\omega^{-1/2}\ee^{-i\pi/4}
                +O(\omega^{-3/2})
\end{equation}\
where \(\chi=\frac b{\sqrt\kappa}.\)
The controller uses the known diffusion factor after separating
$\chi$. Changing the sign of $b$ changes the total transfer phase
by $\pi$, but does not change the phase of this factor.
{The demodulation phase is thus fixed in advance, and $\chi$
is absorbed into a loop gain  shown to be proportional to the positive
quantity $b^2/\sqrt\kappa$ in Section~\ref{sec:stability}.}

%% file: sections/controller.tex
\section{Controller and Averaged Target}\label{sec:controller}
Fix $a_0>0$, a dither amplitude $a\in(0,a_0]$, and $\lambda>0$.
Choose gains satisfying
\begin{equation}\label{eq:gaincondition}
 0<k_I<\lambda k_P.
\end{equation}
Set $\theta=\omega(t-t_0)+\theta_0$ and $\eps=\omega^{-1/2}$.
The two-state controller is
\begin{equation}\label{eq:piarchitecture}
 u=v-k_P\widehat g+a\sin\theta,\qquad \dot v=-k_I\widehat g,
\end{equation}
with the implemented filter
\begin{equation}\label{eq:nominal}
\begin{aligned}
 m&=\frac{2\lambda\eps}{a}\sin(\theta+\pi/4),\qquad
 \dot m=\frac{2\lambda}{a\eps}\cos(\theta+\pi/4),\\
 \widehat g&=p+mJ,\qquad
 \dot p=-\lambda p-(\lambda m+\dot m)J.
\end{aligned}
\end{equation}
Only the cost \eqref{eq:cost}, the states $(v,p)$, and the dither
are used. The known derivative $\dot m$ is computed from time.
For absolutely continuous cost traces, the equations imply
\begin{equation}\label{eq:filteridentity}
 \dot{\widehat g}=-\lambda\widehat g+m\dot J.
\end{equation}
This is an analysis identity, not an implementation requiring
cost differentiation. Section~\ref{sec:regularity} supplies the
needed trace regularity.

{The controller is a proportional--integral extremum
seeking controller. The bias $v$ is the slowly adapted heater power, $-k_P\widehat
g$ is a proportional correction, and $a\sin\theta$ is the probing
signal. The filter \eqref{eq:nominal} estimates the sensitivity of the
cost to the heater power. By \eqref{eq:filteridentity}, it multiplies
the rate of change of the measured cost by a copy of the dither
advanced by $\pi/4$, which is exactly in phase with the rate of change
of the cost ripple produced by the thermal wave. It passes the
product through a first-order low-pass filter of bandwidth $\lambda$.
The factor $\eps/a$ in $m$ normalizes the product. The differentiation of a
ripple of amplitude $a\chi\eps\varphi'$ produces a factor
$a\omega^{1/2}$, that is canceled by $\eps/a$.}

Define the signed sensitivity
\begin{equation}\label{eq:sensitivity}
 g(x)=\chi\varphi'(\Cb x).
\end{equation}
At a frozen boundary value $z$, the leading cost oscillation is
$a\chi\eps\varphi'(z)\sin(\theta-\pi/4)$.
Multiplication of its derivative by $m$ has mean
$\lambda\chi\varphi'(z)$. This explains the modulation in
\eqref{eq:nominal}, and the exact argument follows below.

Using the sensitivity in \eqref{eq:sensitivity}, the averaged target
retains the finite filter bandwidth:
\begin{equation}\label{eq:average}
\begin{aligned}
 \bar x_t&=\kappa\Delta\bar x-\mu_0\bar x+F(s,\bar x)+d\\
 -\kappa\bar x_{s_1}(0,s_2,t)&=bq(s_2)(\bar v-k_P\bar\xi)\\
 \dot{\bar v}&=-k_I\bar\xi\\
 \dot{\bar\xi}&=-\lambda\bar\xi+\lambda g(\bar x).
\end{aligned}
\end{equation}
The other edges have the conditions in \eqref{eq:plant}.
{The averaged target is the loop that would be obtained if the
sensitivity $g$ were measured directly and passed through the same
low-pass filter. The dither and the demodulation have disappeared, and
the task of the following sections is to show how closely the
implemented loop follows it.}
Write $\bar Z=(\bar x,\bar v,\bar\xi)$,
$\Zs=(w^\star,u^\star,0)$, and $X_\omega=(w,v,\widehat g)$.
The common state space $\E=\X\times\R^2$ has norm
\begin{equation}\label{eq:statespace}
 \norm{(e,\nu,\zeta)}_\E=\norm e_\infty+|\nu|+|\zeta|,
\end{equation}
and the regulation error is
\begin{equation}\label{eq:error}
 E_\omega(t)=\norm{X_\omega(t)-\Zs}_\E.
\end{equation}
We seek local exponential transients with ultimate state and cost
bounds of orders $\omega^{-1/2}$ and $\omega^{-1}$, respectively,
at fixed $a$. 

A constant cost offset $J_c=J+c$ leaves the implemented estimate
and physical trajectories unchanged under $p_c=p-mc$, including
$p_c(t_0)=p(t_0)-m(t_0)c$. An offset-invariant initialization is
\begin{equation}\label{eq:filterinitialization}
 p(t_0)=-m(t_0)J(t_0),\qquad \widehat g(t_0)=0.
\end{equation}
Accordingly, \eqref{eq:error} measures $\widehat g$, whereas
$p=\widehat g-mJ$ can depend on the absolute cost level.
{For the optical filter, this means that the controller does not
need to know the peak transmission of the element.}

%% file: sections/averaged_stability.tex
\section{Stability of the Averaged System}\label{sec:stability}

For the averaged state $Z=(x,v,\xi)$, set
$Z_e=(e,\nu,\xi)=(x-w^\star,v-u^\star,\xi)$.
We retain a longitudinal linear contribution to the reaction, prove
stability of the resulting field--controller system, and then bound
the remaining nonlinear terms in the norm \eqref{eq:statespace}.
The linear estimate includes the sensitivity filter and supplies
both a decay rate and a gain for additional boundary commands.
These two estimates distinguish boundary disturbances from
integrable internal forcing in the nonlinear robustness result.

{The argument has three parts. First, a linear model is extracted
that keeps the diffusion, the part of the reaction that helps (through
a reference slope $V_0$), and the controller. Its stability is proved
with an energy function that pairs the mean temperature across the
plate with the two controller states. Second, the response of this
linear model to an additional boundary command is bounded. Third, what
the linear model leaves out is shown to be a small Lipschitz
perturbation, which can be absorbed using a Gronwall-type argument. 
The unknown sign of $b$ enters twice, once in the effect of the heater on the plate and once in the
reconstructed sensitivity. As a result, the loop gain seen by the energy estimate
is $b^2h/\sqrt\kappa>0$.}

\subsection{Linear system and constants}
Choose $V_0\in C([0,1])$ to include stabilizing reaction slopes in
the linear calculation. The remainder is
$F(s,w^\star+e)-F(s,w^\star)-V_0(s_1)e$.
The controller does not use $V_0$. We define
\begin{equation}\label{eq:referencepotential}
\begin{aligned}
 \mu_{\rm ref}(s_1)&=\mu_0-V_0(s_1),\qquad
 \underline\mu=\min_{[0,1]}\mu_{\rm ref},\\
 \mu_{\rm av}&=\int_0^1\mu_{\rm ref}(y)\dd y.
\end{aligned}
\end{equation}
\begin{assumption}\label{ass:damping}
The damping in \eqref{eq:referencepotential} satisfies
$\underline\mu>0$.
\end{assumption}
Linearizing the sensitivity at the optimum gives
\begin{equation}\label{eq:linearreference}
\begin{aligned}
 e_t&=\kappa\Delta e-\mu_{\rm ref}(s_1)e\\
 -\kappa e_{s_1}(0,s_2,t)&=bq(s_2)(\nu-k_P\xi+d_u)\\
 \dot\nu&=-k_I\xi\\
 \dot\xi&=-\lambda\xi+\lambda\chi h\Cb e
\end{aligned}
\end{equation}
with the boundary conditions descric=bed in \eqref{eq:plant}.
Let $T_0$ be its unforced semigroup on $\E$, and let
$S_{\rm ref}$ and $\KBref$ be the homogeneous Neumann semigroup
of $A+V_0$ and its boundary-input kernel.

We set $c_\star=h/\sqrt\kappa$, $c=b^2c_\star$,
$D=k_P-k_I/\lambda>0$, $Q=1+|\vartheta|$, and
$\Jk(\beta)=\beta^{-1}+(\kappa\beta)^{-1/2}$ for $\beta>0$.
The following constant is chosen
\begin{equation}\label{eq:strictification}
\begin{aligned}
0<\delta_L<\min\Bigl\{&\tfrac12\sqrt{c/k_I},\ \kappa/k_I,\\
 &\left(3k_I/\underline\mu+2\mu_{\rm av}/c\right)^{-1},\ (2D)^{-1}\Bigr\}.
\end{aligned}
\end{equation}
We then define the constants
\begin{equation}\label{eq:referenceconstants}
\begin{aligned}
 \sigma&=\tfrac13\min\{\underline\mu,\delta_Lk_I,\lambda\},\\
 L_\perp&=\sqrt{ck_I}+\sqrt{cD\lambda}+\delta_Lk_I,
\end{aligned}
\end{equation}
\begin{equation}\label{eq:controllerbounds}
\begin{aligned}
 L_U&=2\sqrt{c_\star\lambda k_P},\,
 L_\nu=2\sqrt{\frac{c_\star k_Ik_P}{D}},\,
 L_\xi=2\sqrt{\frac{c_\star\lambda}{D}}.
\end{aligned}
\end{equation}
and
\begin{equation}\label{eq:initialenergy}
\begin{aligned}
 P_c&=c_\star^{-1}
 \begin{pmatrix}1/k_I&-1/\lambda\\-1/\lambda&k_P/\lambda\end{pmatrix},\\
 K_W&=\max\{1,\lambda_{\max}(P_c)\},\\
 C_E&=\frac{\sqrt{3K_W}}2+\frac{4L_\perp}{\underline\mu-\sigma}.
\end{aligned}
\end{equation}
Here $P_c>0$ follows from $\det P_c=D/(c_\star^2k_I\lambda)>0$,
and $\sigma\le\underline\mu/3$. Finally, we let
\begin{equation}\label{eq:M}
\begin{aligned}
 M={}&1+C_E\bigl[|b|QL_U\Jk(\underline\mu-\sigma)+L_\nu+L_\xi\bigr],
\end{aligned}
\end{equation}
\[
 H_B=b^2(\kappa^{-1}+\underline\mu^{-1}+\delta_L/c),\quad
 G_V=\sqrt{H_B/\sigma},
\]
and
\begin{equation}\label{eq:GBdef}
\begin{aligned}
 G_B={}&|b|Q\Jk(\underline\mu)(1+L_UG_V)+(L_\nu+L_\xi)G_V.
\end{aligned}
\end{equation}
These constants are independent of the nonlinear radius chosen below.
Their dependence on the plant and the chosen reaction coefficient
is explicit. The restriction on $V_0$ separates the linear
transverse modes. It does not require the full nonlinear reaction
$F$ to be independent of $s_2$.

\begin{proposition}\label{prop:linear}
Under Assumptions~\ref{ass:optimum}, \ref{ass:damping}, and
\eqref{eq:gaincondition}, system \eqref{eq:linearreference}
defines a strongly continuous semigroup satisfying
\begin{equation}\label{eq:linear}
 \norm{T_0(t)}_{\mathcal L(\E)}\le M\ee^{-\sigma t},\qquad t\ge0,
\end{equation}
with \eqref{eq:referenceconstants} and \eqref{eq:M}.
Its zero-state response to $d_u\in L^\infty_{\rm loc}$ satisfies
\begin{equation}\label{eq:GB}
 \norm{Z_{d_u}(t)}_\E\le G_B\norm{d_u}_{\infty,[t_0,t]},
 \qquad t\ge t_0,
\end{equation}
with \eqref{eq:GBdef} and the essential-supremum input norm.
Continuous initial fields require no boundary compatibility.
\end{proposition}

\begin{proof} The proof proceeds in three steps. 
\emph{1. Kernel and transverse modes.}
Let $K_N,K_{\rm ref}$ denote the interval Neumann kernels of
$\kappa\partial_{yy}$ and $\kappa\partial_{yy}-\mu_{\rm ref}$.
For nonnegative data, the Neumann heat solution multiplied by
$\ee^{-\underline\mu t}$ is a supersolution. Comparison gives
$0\le K_{\rm ref}\le\ee^{-\underline\mu t}K_N$ and
$\|S_{\rm ref}(t)\|\le\ee^{-\underline\mu t}$
\cite[Sec.~2.3]{Ouhabaz2005}.
Its residual is $(\mu_{\rm ref}-\underline\mu)$ times the
nonnegative supersolution. Positivity and preservation of
constants by the heat semigroup then give the asserted norm
bound for arbitrary signed initial fields.
Green's identity and separation give
\begin{equation}\label{eq:kernel}
\begin{aligned}
 K_N(t;y,0)&=1+2\sum_{n\ge1}\ee^{-\kappa n^2\pi^2t}\cos(n\pi y)\\
 \KBref(t)(s)&=bK_{\rm ref}(t;s_1,0)[1+\vartheta\ee^{-\kappa\pi^2t}\cos(\pi s_2)]\\
 \norm{\KBref(t)}_\infty
 &\le |b|Q\ee^{-\underline\mu t}(1+1/\sqrt{\pi\kappa t}).
\end{aligned}
\end{equation}
Indeed, the Gaussian sum is at most its integral over
$[0,\infty)$, namely $(2\sqrt{\pi\kappa t})^{-1}$.
Thus, for $\zeta<\underline\mu$, we get
\begin{equation}\label{eq:kernelintegral}
 \int_0^\infty\ee^{\zeta t}\norm{\KBref(t)}_\infty\dd t
 \le |b|Q\Jk(\underline\mu-\zeta).
\end{equation}
For example, the field component of the mild equation is
\[
\begin{aligned}
 e(t)={}&S_{\rm ref}(t-t_0)e(t_0)\\
 &+\int_{t_0}^t\KBref(t-s)
              [\nu(s)-k_P\xi(s)+d_u(s)]\dd s.
\end{aligned}
\]
The controller components are ordinary integral equations.
Their Lipschitz constants tend to zero with the interval length,
as does $\|\KBref\|_{L^1(0,T;\X)}$.
The coupled equations therefore define a Volterra contraction
for small time. The linear coefficients permit a uniform restart
interval, giving global mild solutions. For zero input, uniqueness,
autonomy, and continuity at the initial time give $T_0$.

The boundary profile contains a mean and one cosine mode.
The actuation therefore forces only two transverse components,
whereas an arbitrary initial field can contain every mode.
We separate these components to identify the part seen by the
sensor and retain the contribution of the remaining field.
We write $\tau=t-t_0$ and define
\[
\begin{aligned}
 e_0(s_1)&=\int_0^1e(s_1,s_2)\dd s_2,\\
 e_1(s_1)&=2\int_0^1e(s_1,s_2)\cos(\pi s_2)\dd s_2,\\
 e_\perp&=e-e_0-e_1\cos(\pi s_2).
\end{aligned}
\]
These projections commute with $S_{\rm ref}$ because $V_0$
is longitudinal. Only $e_0,e_1$ are forced. The first cosine
vanishes at the sensor, so \(e_\perp(t)=S_{\rm ref}(\tau)e_\perp(t_0)\), \(z_\perp=\Cb e_\perp\), \( \Cb e=e_0(0)+z_\perp\) and therefore
\begin{equation}\label{eq:unforcedremainder}
\begin{aligned}
 |z_\perp(t)|&\le4\ee^{-\underline\mu\tau}\norm{e(t_0)}_\infty.
\end{aligned}
\end{equation}
This expression bounds the remainder projection on $\X$,
since the identity, mean, and first-cosine projections have
norms at most $1,1,2$. Although $e_1$ is forced, it contributes
nothing to the midpoint observation.
Thus $z_\perp$ acts as an exponentially decaying input to the
mean-field and controller subsystem.

\emph{Energy and full-field decay.}
Initially take smooth compatible data and $d_u=0$. Put
$\eta_c=\xi/(bc_\star)$ and
$\rho_c=-\nu/(bc_\star)+(k_I/\lambda)\eta_c$.
The mean flux is $-c(\rho_c+D\eta_c)$, and
\[
\begin{aligned}
 \dot\rho_c&=k_I[e_0(0)+z_\perp],\\
 \dot\eta_c&=\lambda[e_0(0)+z_\perp-\eta_c].
\end{aligned}
\]
These coordinates use $\chi h=bc_\star$. The same unknown sign
enters the plant flux and the reconstructed sensitivity, so their
product in the feedback dynamics is $c=b^2c_\star>0$.
The coordinate change is used for the stability proof.
We define $j_0=\int_0^1e_0$, $I_e=\int_0^1\mu_{\rm ref}e_0^2$,
and $j_\mu=\int_0^1\mu_{\rm ref}e_0$, where
$\|\cdot\|_2$ the $L^2(0,1)$ norm. We pose
\begin{equation}\label{eq:storage}
\begin{aligned}
 W&=\tfrac12\norm{e_0}_2^2+\frac{c\rho_c^2}{2k_I}
           +\frac{cD\eta_c^2}{2\lambda},\\
 V&=W+\delta_L\rho_cj_0.
\end{aligned}
\end{equation}
Since $|\rho_cj_0|\le\sqrt{k_I/c}\,W$,
\eqref{eq:strictification} gives $W/2\le V\le3W/2$.
The energy $W$ pairs the mean field with the controller coordinates
so that the terms involving $e_0(0)$ cancel under integration by
parts. Its derivative does not directly penalize $\rho_c$.
The additional term $\delta_L\rho_cj_0$ uses the mean balance
to obtain that dissipation while keeping $V$ comparable to $W$.
Integration by parts yields
\[
\begin{aligned}
 \dot W&=-\kappa\norm{e_0'}_2^2-I_e-cD\eta_c^2 +(c\rho_c+cD\eta_c)z_\perp,\\
 \dot j_0&=-c\rho_c-cD\eta_c-j_\mu.
\end{aligned}
\]
Here primes denote spatial derivatives.
The identity \[e_0(0)=j_0-\int_0^1(1-y)e_0'(y)\dd y\] gives
\[|e_0(0)j_0|\le(3\|e_0\|_2^2+\|e_0'\|_2^2)/2.\]
Applying the Cauchy--Schwarz inequality to $\sqrt{\mu_{\rm ref}}$ and
$\sqrt{\mu_{\rm ref}}e_0$ gives $j_\mu^2\le\mu_{\rm av}I_e$
\cite[Sec.~1.1]{HardyLittlewoodPolya1952}.
Differentiating \eqref{eq:storage} and using
$I_e\ge\underline\mu\|e_0\|_2^2$, the mixed terms are shown to satisfy
\[
\begin{aligned}
 \delta_L k_I|e_0(0)j_0|
 &\le\tfrac{\delta_Lk_I}{2}\norm{e_0'}_2^2
                    +\tfrac{3\delta_Lk_I}{2\underline\mu}I_e,\\
 \delta_LcD|\rho_c\eta_c|
 &\le\tfrac{\delta_Lc}{4}\rho_c^2+\delta_LcD^2\eta_c^2,\\
 \delta_L|\rho_cj_\mu|
 &\le\tfrac{\delta_Lc}{4}\rho_c^2+\tfrac{\delta_L\mu_{\rm av}}c I_e.
\end{aligned}
\]
Each inequality accounts for a specified part of the dissipation.
The condition involving $\mu_{\rm av}$ controls the weighted
reaction integral, $\delta_L<\kappa/k_I$ controls the field
derivative term, and $\delta_L<(2D)^{-1}$ controls the mixed
controller term.
The derivative of the cross term contributes
$-\delta_Lc\rho_c^2$. These bounds account for half of this term
and at least half of each original dissipative term by
\eqref{eq:strictification}. Hence
\begin{equation}\label{eq:energy}
\begin{aligned}
 \dot V\le{}&-\tfrac\kappa2\norm{e_0'}_2^2-\tfrac12I_e
           -\tfrac{\delta_Lc}{2}\rho_c^2-\tfrac{cD}{2}\eta_c^2\\
 &+2L_\perp\sqrt V\,|z_\perp|.
\end{aligned}
\end{equation}
Indeed, $V\ge W/2$ bounds the forcing coefficient
$|c\rho_c+cD\eta_c+\delta_Lk_Ij_0|$ by $2L_\perp\sqrt V$.
The negative terms in \eqref{eq:energy} dominate
$3\sigma W\ge2\sigma V$.
The comparison for $\sqrt V$, justified at zero by positive
regularization, gives
\[
\begin{aligned}
 \sqrt{V(t)}\le{}&\ee^{-\sigma\tau}\sqrt{V(t_0)}+L_\perp\int_{t_0}^t
                 \ee^{-\sigma(t-s)}|z_\perp(s)|\dd s.
\end{aligned}
\]

The energy estimate controls the mean field in $L^2$ and the
controller states. To recover the product norm, we first bound
the feedback command $\nu-k_P\xi$. Its action on the full field
is then estimated through the boundary kernel, including both
the forced first cosine and the remaining modes.
For arbitrary real $a_1,a_2$, the controller terms in
$V\ge W/2$ give
\[
 |a_1\rho_c+a_2\eta_c|^2
 \le4\left(\frac{k_Ia_1^2}c+\frac{\lambda a_2^2}{cD}\right)V.
\]
We then use $\nu-k_P\xi=-bc_\star(\rho_c+D\eta_c)$,
$\nu=bc_\star[-\rho_c+(k_I/\lambda)\eta_c]$, and
$\xi=bc_\star\eta_c$ to generate the bounds 
\[
 |\nu-k_P\xi|\le L_U\sqrt V,\quad
 |\nu|\le L_\nu\sqrt V,\quad |\xi|\le L_\xi\sqrt V,
\]
with \eqref{eq:controllerbounds}. Moreover,
$W=\|e_0\|_2^2/2+(\nu,\xi)P_c(\nu,\xi)^\top/2$.
Thus, for $N_0=\|Z_e(t_0)\|_\E$,
$\sqrt{V(t_0)}\le\sqrt{3K_W}N_0/2$.
Using \eqref{eq:unforcedremainder} in scalar comparison gives
$\sqrt{V(t)}\le C_E\ee^{-\sigma\tau}N_0$, with
\eqref{eq:initialenergy}, since
$\int_0^\tau\ee^{-(\underline\mu-\sigma)s}\dd s
\le(\underline\mu-\sigma)^{-1}$.
The two terms in $C_E$ are as follows. The first term
converts the initial product norm into the energy norm.
The second term accounts for the decaying observation of the
unforced transverse modes. Its denominator is positive
since these modes are estimated with the faster heat-semigroup
decay $\underline\mu>\sigma$.
The field mild equation and
\eqref{eq:kernelintegral} imply
\[
\begin{aligned}
 \norm{e(t)}_\infty
 &\le \ee^{-\sigma\tau}N_0+|b|QL_UC_E\Jk(\underline\mu-\sigma)
                          \ee^{-\sigma\tau}N_0.
\end{aligned}
\]
Adding the controller bounds proves \eqref{eq:linear}.

\emph{2. Boundary gain and mild data.}

For zero initial state, we have $z_\perp=0$. The additional boundary
flux makes $\dot V$ gain $b[e_0(0)+\delta_L\rho_c]d_u$.
Using $|e_0(0)|\le\|e_0\|_2+\|e_0'\|_2$, Young's inequality gives
\[
\begin{aligned}
 |b|\norm{e_0'}_2|d_u|
 &\le\tfrac\kappa4\norm{e_0'}_2^2+\tfrac{b^2}\kappa|d_u|^2,\\
 |b|\norm{e_0}_2|d_u|
 &\le\tfrac14I_e+\tfrac{b^2}{\underline\mu}|d_u|^2,\\
 |b|\delta_L|\rho_c||d_u|
 &\le\tfrac{\delta_Lc}{4}\rho_c^2+\tfrac{b^2\delta_L}{c}|d_u|^2.
\end{aligned}
\]
The remaining dissipation dominates $3\sigma W/2\ge\sigma V$,
so $\dot V\le-\sigma V+H_B|d_u|^2$ and
$\sqrt V\le G_V\|d_u\|_{\infty,[t_0,t]}$.
Here integration from $V(t_0)=0$ bounds $V(t)$ by
$H_B\|d_u\|_{\infty,[t_0,t]}^2/\sigma$.
Applying \eqref{eq:kernelintegral} to the total command
$\nu-k_P\xi+d_u$ yields
\[
\begin{aligned}
 \norm{e(t)}_\infty
 &\le |b|Q\Jk(\underline\mu)(1+L_UG_V)\norm{d_u}_{\infty,[t_0,t]}.
\end{aligned}
\]
The controller coordinates contribute at most
\[
 (L_\nu+L_\xi)G_V\|d_u\|_{\infty,[t_0,t]}.
\]
Adding the field and controller bounds gives \eqref{eq:GB}
with exactly \eqref{eq:GBdef}.

{\emph{3. Passage to mild data.}}
The integration by parts above uses smooth compatible data, but
the resulting constants do not depend on higher derivatives of
those data. We approximate the initial field for the unforced
semigroup estimate and the input for the zero-state boundary
estimate. First we approximate a continuous initial field uniformly
using tensor cosine sums. We add
$-(b/\kappa)U_0q(s_2)s_1\eta(ns_1)$ to impose the initial flux,
where $U_0=\nu(t_0)-k_P\xi(t_0)$ and
$\eta\in C^\infty([0,\infty))$ equals one near zero and vanishes
on $[1,\infty)$. This correction tends uniformly to zero and
preserves the other Neumann conditions since $q'(0)=q'(1)=0$.
Its derivative at the actuated edge is exactly
$-(b/\kappa)U_0q(s_2)$, while its supremum norm is $O(n^{-1})$.
On each finite interval, approximate bounded inputs in $L^p$,
$p>2$, by smooth inputs with no larger supremum norm, vanishing
near $t_0$ for the zero-state estimate. Since $p'=p/(p-1)<2$,
\eqref{eq:kernel} gives $\KBref\in L^{p'}(0,T;\X)$.
The vanishing near $t_0$ makes the zero initial field compatible
before the approximating flux begins. Upon preserving the input norm,
the right-hand side of the boundary-gain estimate remains bounded
by the norm of the original input.
Using H\"older's inequality, the convolution difference is bounded by
$\|\KBref\|_{L^{p'}}\|d_{u,n}-d_u\|_{L^p}\to0$, uniformly in
time. The Volterra continuous dependence leads 
\eqref{eq:linear} and \eqref{eq:GB} to mild solutions.
Compatibility is imposed only on the approximating data and does
not restrict the initial fields in the proposition.
\end{proof}

\subsection{Nonlinear margin and robustness}
\begin{assumption}\label{ass:margin}
There exist $R,\eta_\star>0$ such that
Assumption~\ref{ass:regularity} holds for fields with
$\|x-w^\star\|_\infty\le R+\eta_\star$, the cost is $C^2$ on
an open interval containing
$[z^\star-R-\eta_\star,z^\star+R+\eta_\star]$, and
\begin{equation}\label{eq:ell}
\begin{aligned}
 \ell(R):={}&\sup_{\substack{s\in\overline\Omega\\|\eta|\le R}}
 |F_w(s,w^\star(s)+\eta)-V_0(s_1)|\\
 &+\lambda|\chi|\sup_{|\eta|\le R}
        |\varphi''(z^\star+\eta)-h|<\sigma/M.
\end{aligned}
\end{equation}
\end{assumption}
Define the nonlinear decay rate and cost-derivative bounds by
\[
\begin{aligned}
 \alpha&=\sigma-M\ell(R)>0,\\
 \Phi_j&=\max_{|z-z^\star|\le R+\eta_\star}|\varphi^{(j)}(z)|,
 \qquad j=1,2.
\end{aligned}
\]
If $F_w(s,w^\star(s))=V_0(s_1)$, continuity gives
$\ell(R)\to0$ as $R\downarrow0$. Otherwise the discrepancy,
including transverse dependence, must satisfy \eqref{eq:ell}.
The larger radius accommodates the harmonic lifting later.
The margin controls the reaction left after $V_0$ is retained,
so large stabilizing equilibrium slopes need not be small.
Plant dissipativity alone does not establish this condition
for the coupled feedback system.
All quantities in the right-hand side of \eqref{eq:ell}
have already been fixed by the linear calculation. Choosing
$R$ therefore introduces no circular dependence into $M$
or $\sigma$. The bounds $\Phi_1,\Phi_2$ are finite by
continuity on the enlarged compact trace interval.

\begin{lemma}\label{lem:averaged}
Under Assumptions~\ref{ass:regularity}--\ref{ass:margin} and
\eqref{eq:gaincondition}, let $\widetilde Z$ solve
\eqref{eq:average}. Let $Z$ solve it with additional terms
$d_w,d_v,d_\xi$ in the field, bias, and filter equations and
boundary command $v-k_P\xi+d_u$. Assume strongly measurable
disturbances with
$ q_d=\|d_w\|_\infty+|d_v|+|d_\xi|\in L^1_{\rm loc}$ and
$d_u\in L^\infty_{\rm loc}$. Writing
$D_Z(t)=\|Z(t)-\widetilde Z(t)\|_\E$, while both states remain
in $B_R(\Zs)$,
\begin{equation}\label{eq:master}
\begin{aligned}
 D_Z(t)\le{}&M\ee^{-\alpha(t-t_0)}D_Z(t_0)
       +G_B\frac\sigma\alpha\norm{d_u}_{\infty,[t_0,t]}\\
 &+M\int_{t_0}^t\ee^{-\alpha(t-s)}q_d(s)\dd s.
\end{aligned}
\end{equation}
\end{lemma}
\begin{proof}
The nonlinear remainder relative to \eqref{eq:linearreference} is
\[
 \mathcal N(Z_e)=
 \begin{pmatrix}
 F(w^\star+e)-F(w^\star)-V_0e\\
 0\\
 \lambda\chi[\varphi'(z^\star+\Cb e)-h\Cb e]
 \end{pmatrix}.
\]
The mean-value formula and \eqref{eq:ell} give
$\mathcal N(0)=0$ and $\operatorname{Lip}(\mathcal N)\le\ell(R)$
on the convex ball $B_R(0)$. Every mean-value segment stays
within the stated ranges, and $\|\Cb\|=1$ bounds the scalar
part by the field difference. Variation of constants relative
to $T_0$ gives local well-posedness by contraction. The boundary
response is continuous and the internal forcing is locally
integrable. The same argument permits continuation while states
remain bounded and strictly inside the regularity region.
Subtracting the mild equations and using
\eqref{eq:linear} and \eqref{eq:GB} yields
\[
\begin{aligned}
 D_Z(t)\le{}&M\ee^{-\sigma(t-t_0)}D_Z(t_0)
                       +G_B\norm{d_u}_{\infty,[t_0,t]}\\
 &+M\int_{t_0}^t\ee^{-\sigma(t-s)}
                  [\ell(R)D_Z(s)+q_d(s)]\dd s.
\end{aligned}
\]
The kernel $k_0(t)=M\ell(R)\ee^{-\sigma t}$ has convolution
powers and resolvent
\[
\begin{aligned}
 k_0^{*n}(t)&=\frac{[M\ell(R)]^n}{(n-1)!}t^{n-1}\ee^{-\sigma t},\\
 r_0(t)&=\sum_{n\ge1}k_0^{*n}(t)=M\ell(R)\ee^{-\alpha t}.
\end{aligned}
\]
Iteration of the inequality therefore adds $r_0$ convolved
with its free terms. The remainder converges to zero on every
finite interval because of the factorial denominator. The identity
$\ee^{-\sigma t}+(r_0*\ee^{-\sigma\cdot})(t)=\ee^{-\alpha t}$
gives the initial and internal terms in \eqref{eq:master}.
The boundary running supremum is nondecreasing, so its factor
is at most $G_B[1+M\ell(R)/\alpha]=G_B\sigma/\alpha$.
Indeed, its convolution with $r_0$ is bounded by its value
at the current time multiplied by
$\int_0^{t-t_0}r_0(s)\dd s\le M\ell(R)/\alpha$.
This also covers $\ell(R)=0$, when $r_0=0$.
\end{proof}

\begin{corollary}\label{cor:target}
Under Lemma~\ref{lem:averaged}, every averaged initial state with
$\|\bar Z(t_0)-\Zs\|_\E<R/M$ generates a unique global mild
solution satisfying
\begin{equation}\label{eq:targetbound}
\begin{aligned}
 \norm{\bar Z(t)-\Zs}_\E
 &\le M\ee^{-\alpha(t-t_0)}\norm{\bar Z(t_0)-\Zs}_\E,\\
 &\hspace{10em}t\ge t_0.
\end{aligned}
\end{equation}
Only continuity of the initial field is required.
\end{corollary}
\begin{proof}
Set $\widetilde Z=\Zs$ and all disturbances to zero in
\eqref{eq:master}. The resulting \eqref{eq:targetbound} stays
strictly below $R$, excluding a first exit by continuity.
The strict initial inequality gives the fixed gap
$R-M\|\bar Z(t_0)-\Zs\|_\E>0$.
This gap persists on every interval of existence, so the
trajectory cannot approach the boundary of the certified
region while remaining bounded inside it.
On a finite maximal interval, the initial-state term in the
mild equation has an endpoint limit by strong continuity of
$T_0$. The nonlinear convolution also has a limit because
$\mathcal N$ is bounded along the trajectory and dominated
convergence applies to the semigroup factor.
Both the field and controller coordinates therefore have
limits in the product space. This limit remains
strictly inside the regularity region and restarts the local
solution, excluding finite-time loss of existence.
\end{proof}

Corollary~\ref{cor:target} gives the decaying averaged
trajectory used below. Lemma~\ref{lem:averaged} controls
the effect of the boundary and internal residuals produced
by the exact transformation of the implemented system.
{The next three sections construct that transformation and bound
its residuals. Section~\ref{sec:harmonic} removes the dither from the
boundary condition, Section~\ref{sec:regularity} establishes the trace
regularity that the filter identity \eqref{eq:filteridentity}
requires, and Section~\ref{sec:comparison} writes the transformed
implemented loop as the averaged target plus explicit disturbances.}

%% file: sections/harmonic_lifting.tex
\section{Harmonic Lifting and Coordinates}\label{sec:harmonic}
{Physically, the dither launches a thermal wave into the plate
whose shape is known exactly for the linear part of the model. There is a
damped oscillation that penetrates a distance of order
$\sqrt{\kappa/\omega}$ from the heated edge and whose amplitude at the
edge is of order $a\omega^{-1/2}$. Subtracting this wave from the
temperature field leaves only the slowly varying part that the
averaged target describes.}
The fixed-amplitude boundary dither produces a field whose supremum
norm decreases with frequency. Then the subtraction of its exact linear response
removes the imposed oscillation from the boundary condition.
For $q_0=1$, $q_1=\vartheta$, we define
\begin{equation}\label{eq:harmonicfield}
\begin{aligned}
 \Psi_\omega(s)&=(i\omega-A_{-1})^{-1}B\\
 &=\frac b\kappa\sum_{j=0}^1q_j\cos(j\pi s_2)
 \frac{\cosh(\gamma_j(1-s_1))}{\gamma_j\sinh\gamma_j}\\
 \gamma_j&=\sqrt{(\mu_0+i\omega)/\kappa+j^2\pi^2}.
\end{aligned}
\end{equation}
All roots have positive real part. The field
$r_\omega=a\operatorname{Im}(\ee^{i\theta}\Psi_\omega)$ solves the
linear diffusion equation with boundary command $a\sin\theta$.
Thus its mild equation subtracts exactly from \eqref{eq:mild}.
The lifting uses the original $A$.

Set $c_h=\coth1$, $\omega_{\rm lift}=\max\{1,2\kappa\}$, and \( L_\rho=|\chi|c_h\),
\begin{equation}\label{eq:liftconstants}
\begin{aligned}
 L_r&=|b|Qc_h/\sqrt\kappa\\
 D_r&=|b|c_h\bigl(Q/\kappa+\pi|\vartheta|/\sqrt\kappa\bigr),\\
 L_{r,\gamma}&=L_r+(2L_r)^{1-\gamma}D_r^\gamma,\\
 L_e&=|\chi|\left(\mu_0/2+4\kappa/(\ee^2-1)\right).
\end{aligned}
\end{equation}
\begin{lemma}\label{lem:lift}
For $\omega\ge\omega_{\rm lift}$, the following holds
\begin{equation}\label{eq:liftbounds}
\begin{aligned}
 \norm{r_\omega}_\infty&\le L_ra\eps,\\
 \norm{r_\omega}_{C^\gamma}&\le L_{r,\gamma}a\eps^{1-\gamma}
 \end{aligned}
\end{equation}
 uniformly in time and phase.
The midpoint trace $\rho_\omega=\Cb r_\omega$ satisfies
\begin{equation}\label{eq:traceharmonic}
\begin{aligned}
 \rho_\omega&=a\chi\eps\sin(\theta-\pi/4)+e_\omega,\\
 |e_\omega|&\le L_ea\eps^3,
\end{aligned}
\end{equation}
and
\begin{equation}\label{eq:rhobounds}
 |\rho_\omega|\le L_\rho a\eps,\qquad
 |\dot\rho_\omega|\le L_\rho a/\eps.
\end{equation}
\end{lemma}
\begin{proof}
Separation of variables gives \eqref{eq:harmonicfield}. The root
bounds $|\gamma_j|\ge\sqrt{\omega/\kappa}$ and
$\operatorname{Re}\gamma_j\ge\sqrt{\omega/(2\kappa)}\ge1$ imply,
for $0\le y\le1$,
\[
 \left|\frac{\cosh(\gamma_j(1-y))}{\sinh\gamma_j}\right|\le c_h,
 \qquad
 \left|\frac{\sinh(\gamma_j(1-y))}{\sinh\gamma_j}\right|\le c_h.
\]
Consequently $\norm{r_\omega}_\infty\le L_ra\eps$ and
$\norm{\nabla r_\omega}_\infty\le D_ra$. The convex-domain
inequality
$[f]_\gamma\le(2\norm f_\infty)^{1-\gamma}
\norm{\nabla f}_\infty^\gamma$, together with $\eps\le1$,
proves \eqref{eq:liftbounds}.

Midpoint evaluation gives \eqref{eq:transfer}, so
$|H_b(i\omega)|\le L_\rho\eps$, which proves
\eqref{eq:rhobounds}. Split the harmonic remainder as
\[
\begin{aligned}
 H_b-\chi(i\omega)^{-1/2}
 ={}&\chi[(\mu_0+i\omega)^{-1/2}-(i\omega)^{-1/2}]\\
 &+\chi(\mu_0+i\omega)^{-1/2}(\coth\gamma_0-1).
\end{aligned}
\]
Integration in $\mu_0$ bounds the first term by
$|\chi|\mu_0\omega^{-3/2}/2$. For
$x_\omega=\sqrt{\omega/(2\kappa)}\ge1$,
\[
 |\coth\gamma_0-1|\le\frac{2\ee^{-2x_\omega}}{1-\ee^{-2}},
 \qquad x_\omega^2\ee^{-2x_\omega}\le\ee^{-2}.
\]
These bounds give the value of $L_e$ in \eqref{eq:liftconstants},
proving \eqref{eq:traceharmonic} and \eqref{eq:asymptotic}.
\end{proof}

{The second coordinate change acts on the filter state. The
product $mJ$ in $\widehat g$ contains a term at twice the dither
frequency and the response to the ripple itself. Both are removed
exactly, so that the remaining variable $\xi$ obeys the averaged filter
equation up to residuals that Section~\ref{sec:comparison} bounds.}

We introduce the triangular change of coordinates
\begin{equation}\label{eq:coordinates}
\begin{aligned}
 x&=w-r_\omega,\\
 z&=\Cb x,\\
 d_{2,\omega}&=\tfrac12\lambda\chi\eps^2\varphi'(z)\cos(2\theta),\\
 \xi&=p+m\varphi(z)-d_{2,\omega},\\
 Z_\omega&=(x,v,\xi).
\end{aligned}
\end{equation}
Its inverse is $w=x+r_\omega$ and
$p=\xi+d_{2,\omega}-m\varphi(\Cb x)$. We define
\[
 \Delta_\omega=d_{2,\omega}
       +m[\varphi(z+\rho_\omega)-\varphi(z)].
\]
Then $\widehat g=\xi+\Delta_\omega$ exactly.

\begin{lemma}\label{lem:coordinates}
Suppose $\omega\ge\omega_{\rm lift}$,
$\norm{x-w^\star}_\infty\le R$, and
$\norm{r_\omega}_\infty\le\eta_\star$. Then
\begin{equation}\label{eq:Deltabound}
 |\Delta_\omega|\le K_\Delta\eps^2
\end{equation}
where \( K_\Delta=\lambda\Phi_1(|\chi|/2+2L_\rho)\). With
\begin{equation}\label{eq:delta}
 \domega=L_ra\eps+K_\Delta\eps^2,
\end{equation}
the coordinate error satisfies
\begin{equation}\label{eq:coordbound}
 \norm{X_\omega-Z_\omega}_\E\le\domega.
\end{equation}
\end{lemma}
\begin{proof}
The two terms in $\Delta_\omega$ are bounded by
$\lambda|\chi|\Phi_1\eps^2/2$ and
$(2\lambda\eps/a)\Phi_1L_\rho a\eps$, respectively.
Adding their sum to the field bound proves
\eqref{eq:coordbound}. Only cost derivatives enter these bounds.
Under a cost offset, $p_c=p-mc$ also leaves $\xi$ unchanged.
\end{proof}

%% file: sections/regularity.tex
\section{Spatial and Boundary-Trace Regularity}\label{sec:regularity}

{This section addresses the filter
identity \eqref{eq:filteridentity}. The analysis multiplies the rate of
change of the measured cost by the demodulating signal, so the time
derivative of the sensor trace $z=\Cb x$ must exist in an integrable
sense. Since the heater acts through the boundary, the response to a
change of command is singular at startup and a bounded derivative
cannot be expected. In what follows, it is shown that the trace is absolutely
continuous and that its derivative, weighted by the averaged decay
rate, is bounded uniformly in time and frequency
(Lemma~\ref{lem:trace}). Readers willing to accept
Lemmas~\ref{lem:kernel}--\ref{lem:trace} can proceed to
Section~\ref{sec:comparison}.}

We establish the spatial bounds and integrable trace derivative needed
for the filter transformation, without initial boundary compatibility.
The boundary kernel is singular at zero, so a bounded pointwise time
derivative is not expected at startup. Spatial H\"older regularity
provides the integrable estimate used below.
Convolutions start at $t_0$, so
$(S*H)(t)=\int_{t_0}^tS(t-s)H(s)\,\dd s$.
With $e=x-w^\star$, $H=F(x+r_\omega)-F(w^\star)$, and
$U=v-u^\star-k_P\widehat g$, the centered mild equation is
\begin{equation}\label{eq:reg_centered_mild}
 e(t)=S(t-t_0)e(t_0)+Y_H(t)+(\KB*U)(t),
\end{equation}
where $Y_H=S*H$.

\subsection{Kernel and spatial estimates}

\begin{lemma}\label{lem:kernel}
For $0<\gamma<1$, define
\begin{equation}\label{eq:reg_kernel_constants}
 \begin{aligned}
 c_S&=2^{1-\gamma}\kappa^{-\gamma/2},\\
 c_A&=2^\gamma(2+\gamma/2)\Gamma(1+\gamma/2),\\
 c_B&=|b|Q(1+2^{\gamma/2}c_S)
                 (1+\sqrt{2/(\pi\kappa)}).
 \end{aligned}
\end{equation}
Here $\Gamma$ is Euler's gamma function.
For $t>0$ and $f\in\X$,
\begin{equation}\label{eq:reg_semigroup_smoothing}
 \norm{S(t)f}_{C^\gamma}
 \le\ee^{-\mu_0t}(1+c_St^{-\gamma/2})\norm f_\infty.
\end{equation}
For $f\in C^\gamma(\overline\Omega)$,
\begin{equation}\label{eq:reg_holder_preservation}
 \norm{S(t)f}_{C^\gamma}
 \le\ee^{-\mu_0t}\norm f_{C^\gamma}.
\end{equation}
The boundary kernel satisfies
\begin{equation}\label{eq:reg_boundary_smoothing}
 \norm{\KB(t)}_{C^\gamma}
 \le c_B\ee^{-\mu_0t}(1+t^{-(1+\gamma)/2}),
\end{equation}
and, for $f\in C^\gamma(\overline\Omega)$,
\begin{equation}\label{eq:reg_generator_smoothing}
 \begin{aligned}
 \norm{AS(t)f}_\infty
 &\le\ee^{-\mu_0t}\bigl[\mu_0\norm f_\infty+c_A\kappa^{\gamma/2}t^{\gamma/2-1}[f]_\gamma\bigr].
 \end{aligned}
\end{equation}
\end{lemma}

\begin{proof}
Extend $f$ evenly and $2$-periodically in each coordinate to
$\widetilde f$ on $\R^2$. The folding map is $1$-Lipschitz, so
$\norm{\widetilde f}_\infty=\norm f_\infty$ and
$[\widetilde f]_\gamma=[f]_\gamma$. With
$\mathcal G_t(r)=(4\pi\kappa t)^{-1}
\exp(-|r|^2/(4\kappa t))$,
\begin{equation}\label{eq:reg_gaussian_representation}
 S(t)f(s)=\ee^{-\mu_0t}
       \int_{\R^2}\mathcal G_t(r)\widetilde f(s-r)\,\dd r.
\end{equation}
The Gaussian has unit mass and
$\norm{\nabla\mathcal G_t}_1=\sqrt\pi/(2\sqrt{\kappa t})
\le(\kappa t)^{-1/2}$.
Convolution therefore bounds the supremum by
$\ee^{-\mu_0t}\norm f_\infty$ and the gradient by
$\ee^{-\mu_0t}(\kappa t)^{-1/2}\norm f_\infty$.
The interpolation inequality
$[g]_\gamma\le(2\norm g_\infty)^{1-\gamma}
\norm{\nabla g}_\infty^\gamma$ on the convex square proves
\eqref{eq:reg_semigroup_smoothing}. Applying
\eqref{eq:reg_gaussian_representation} to differences of translates
proves \eqref{eq:reg_holder_preservation}.
These are estimates on the H\"older norm. They do not require
strong continuity of the semigroup on the full H\"older space.

For the interval Neumann kernel $K_N$, separation gives
\begin{equation}\label{eq:reg_boundary_kernel}
 \begin{aligned}
 \KB(t)(s)&=b\ee^{-\mu_0t}K_N(t,s_1,0)[1+\vartheta\ee^{-\kappa\pi^2t}\cos(\pi s_2)].
 \end{aligned}
\end{equation}
By expansion of the the cosine, its positivity  implies
\[
 0\le K_N(t,y,0)
 \le1+2\int_0^\infty\ee^{-\kappa\pi^2tr^2}\,\dd r
 =1+(\pi\kappa t)^{-1/2}.
\]
Using this bound in \eqref{eq:reg_boundary_kernel}, we then apply
\eqref{eq:reg_semigroup_smoothing} to
$\KB(t)=S(t/2)\KB(t/2)$. The expansion of the factors
$1+2^{\gamma/2}c_St^{-\gamma/2}$ and
$1+\sqrt{2/(\pi\kappa)}t^{-1/2}$ proves
\eqref{eq:reg_boundary_smoothing} with \eqref{eq:reg_kernel_constants}.
Indeed, every power in the product is bounded by
$1+t^{-(1+\gamma)/2}$, and $c_B$ is the sum of its coefficients.

We then write $S_N(t)=\ee^{\mu_0t}S(t)$.
The cancellation of $\int\partial_t\mathcal G_t=0$ gives
\[
 \kappa\Delta S_N(t)f(s)
 =\int_{\R^2}\partial_t\mathcal G_t(r)
       [\widetilde f(s-r)-\widetilde f(s)]\,\dd r.
\]
The Gaussian identities are such that
\[
 \begin{aligned}
 |\partial_t\mathcal G_t(r)|
 &\le[t^{-1}+|r|^2/(4\kappa t^2)]\mathcal G_t(r),\\
 \int_{\R^2}|r|^q\mathcal G_t(r)\,\dd r
 &=(4\kappa t)^{q/2}\Gamma(1+q/2)
 \end{aligned}
\]
with $q=\gamma$ and $q=\gamma+2$ yield
$\int|\partial_t\mathcal G_t(r)|\,|r|^\gamma\,\dd r
\le c_A\kappa^{\gamma/2}t^{\gamma/2-1}$.
Since $AS(t)=\ee^{-\mu_0t}(\kappa\Delta-\mu_0)S_N(t)$,
we prove \eqref{eq:reg_generator_smoothing}.
\end{proof}

The integration of the bounds in Lemma~\ref{lem:kernel} yields
\begin{equation}\label{eq:reg_integrated_constants}
 \begin{aligned}
 I_S&=\mu_0^{-1}
       +c_S\Gamma(1-\gamma/2)\mu_0^{\gamma/2-1},\\
 I_B&=c_B[\mu_0^{-1}
       +\Gamma((1-\gamma)/2)\mu_0^{(\gamma-1)/2}],\\
 I_A&=1+c_A\kappa^{\gamma/2}\Gamma(\gamma/2)\mu_0^{-\gamma/2}.
 \end{aligned}
\end{equation}
They bound the integrals of $\norm{S(t)}_{\X\to C^\gamma}$,
$\norm{\KB(t)}_{C^\gamma}$, and
$\norm{AS(t)}_{C^\gamma\to\X}$, respectively.
All integrals are finite. Since $0<\gamma<1$, the singularities are
integrable at zero, and $\mu_0>0$ controls the tail.

\begin{lemma}\label{lem:spatial}
Suppose Assumptions~\ref{ass:regularity}--\ref{ass:margin} and
\eqref{eq:initialfamily} hold. Fix $0<a\le a_0$ and
$\omega\ge\omega_{\rm lift}\ge1$. Consider a nominal mild
solution on $I_T=[t_0,t_0+T]$ satisfying
\begin{equation}\label{eq:reg_stopped_region}
 Z_\omega(t)\in B_R(\Zs),\qquad
 \norm{r_\omega(t)}_\infty\le\eta_\star.
\end{equation}
Set $G=R+K_\Delta$ and $U_\star=R+k_PG$. Define
\begin{equation}\label{eq:reg_equilibrium_constants}
 \begin{aligned}
 W_\gamma&=I_S\norm{F(w^\star)+d}_\infty+I_B|u^\star|,\\
 E_{\rm in}&=H_0+L_{r,\gamma}a_0,
 \end{aligned}
\end{equation}
and
\begin{equation}\label{eq:reg_spatial_constants}
 \begin{aligned}
 H_\infty&=L_F(R+\eta_\star),\\
 E_\gamma&=E_{\rm in}+I_SH_\infty+I_BU_\star,\\
 H_\gamma&=H_\infty+2H_F+L_F(2W_\gamma+E_\gamma+L_{r,\gamma}a_0).
 \end{aligned}
\end{equation}
Then
\begin{equation}\label{eq:reg_spatial_bounds}
 \begin{aligned}
 \norm{w^\star}_{C^\gamma}&\le W_\gamma,\\
 \sup_{t\in I_T}\norm{e(t)}_{C^\gamma}&\le E_\gamma,\\
 \sup_{t\in I_T}\norm{H(t)}_{C^\gamma}&\le H_\gamma.
 \end{aligned}
\end{equation}
Moreover, $Y_H\in C^1(I_T,\X)$, with a right derivative at $t_0$,
and
\begin{equation}\label{eq:reg_forcing_derivative}
 Y_H'(t)=H(t)+\int_{t_0}^t AS(t-s)H(s)\,\dd s.
\end{equation}
Then 
\begin{equation}\label{eq:reg_time_bounds}
 \begin{aligned}
 \int_0^\infty\norm{AS(s)e(t_0)}_\infty\,\dd s&\le H_{\rm in},\\
 \sup_{t\in I_T}\norm{Y_H'(t)}_\infty&\le L_H.
 \end{aligned}
\end{equation}
for frequency-independent constants
\begin{equation}\label{eq:reg_time_constants}
\begin{aligned}
 H_{\rm in}&=I_AE_{\rm in},\\
 L_H&=H_\infty+I_AH_\gamma.
\end{aligned}
\end{equation}
The constants are independent of $T$, initial time, and phase,
with $a_0$ and the initial-family bound fixed.
\end{lemma}

\begin{proof}
Equation~\eqref{eq:Deltabound} and $\eps\le1$ give
$|\widehat g|\le G$ and $|U|\le U_\star$. The reaction Lipschitz
bound gives $\norm H_\infty\le H_\infty$. Applying Lemma~\ref{lem:kernel}
to the stationary mild equation at time $T_1>0$ yields
\[
 \norm{w^\star}_{C^\gamma}
 \le\ee^{-\mu_0T_1}(1+c_ST_1^{-\gamma/2})
                   \norm{w^\star}_\infty+W_\gamma.
\]
Letting $T_1\to\infty$ proves the equilibrium bound in
\eqref{eq:reg_spatial_bounds}, with \eqref{eq:reg_equilibrium_constants}.
Only boundedness of
$F(w^\star)+d$ is needed.
The initial-family and lift bounds give
$\norm{e(t_0)}_{C^\gamma}\le H_0+L_{r,\gamma}a\eps^{1-\gamma}
\le E_{\rm in}$.
Equations~\eqref{eq:reg_centered_mild},
\eqref{eq:reg_holder_preservation}, and
\eqref{eq:reg_integrated_constants} then give the bound $E_\gamma$.
Finally, \eqref{eq:regularity} implies
\[
 [H]_\gamma\le2H_F+
 L_F(2[w^\star]_\gamma+[e]_\gamma+[r_\omega]_\gamma),
\]
which proves the $H_\gamma$ bound in
\eqref{eq:reg_spatial_constants}. Thus the spatial regularity of
$H$ follows after that of $e$, without a circular bootstrap.

The forcing $H$ is continuous into $\X$. By
\eqref{eq:reg_generator_smoothing},
\[
 \norm{AS(r)H(s)}_\infty
 \le\ee^{-\mu_0r}(\mu_0+c_A\kappa^{\gamma/2}r^{\gamma/2-1})H_\gamma.
\]
This majorant is integrable, with integral $I_AH_\gamma$.
The integrand is strongly measurable into $\X$ for $r>0$.
For each fixed positive $r$, continuity of $H$ into $\X$
suffices because $AS(r)$ is bounded on that space. No continuity
of $H$ in the stronger H\"older norm is assumed.
The time derivative must be justified at the upper integration
limit, where $AS(t-s)$ is singular. The factor $r^{\gamma/2-1}$
is integrable because $\gamma>0$. We therefore integrate the
semigroup derivative first and interchange the resulting integrals.
This produces an integral of a continuous $\X$-valued function,
to which the fundamental theorem of calculus applies.
The identity $\int_0^\tau AS(r)f\,\dd r=S(\tau)f-f$ for
$f\in C^\gamma$ and Fubini's theorem give
\[
 Y_H(t)=\int_{t_0}^t[H(s)+(AS*H)(s)]\,\dd s.
\]
The domination of the convergence with the same majorant makes the integrand
continuous, including its right limit at $t_0$. This proves
\eqref{eq:reg_forcing_derivative} and, by integration of
\eqref{eq:reg_generator_smoothing}, the constants and bounds in
\eqref{eq:reg_time_constants}--\eqref{eq:reg_time_bounds}.
\end{proof}

\subsection{Absolute continuity and the weighted trace bound}

Evaluating \eqref{eq:reg_centered_mild} gives
\begin{equation}\label{eq:tracerepresentation}
 \begin{aligned}
 z&=z_{\rm f}+k*U,\\
 z_{\rm f}&=z^\star+\Cb S(t-t_0)e(t_0)+\Cb Y_H,
 \end{aligned}
\end{equation}
where $k=\Cb\KB=b\ee^{-\mu_0t}K_N(t,0,0)$ has constant sign.
Its norm is
\begin{equation}\label{eq:reg_K0}
 K_0=\norm k_{L^1(0,\infty)}
 =\frac{|b|}{\sqrt{\kappa\mu_0}}\coth\sqrt{\mu_0/\kappa}.
\end{equation}
Indeed, the integral of $k$ is the boundary value of the
stationary solution of $\kappa y''-\mu_0y=0$ with
$-\kappa y'(0)=b$ and $y'(1)=0$, which gives \eqref{eq:reg_K0}.

\begin{lemma}\label{lem:trace}
Under Lemma~\ref{lem:spatial}, $z=\Cb x$ is absolutely continuous
on $I_T$. Define
\begin{equation}\label{eq:reg_trace_coefficients}
 \begin{aligned}
 C_0&=|k_P\lambda-k_I|G+2k_P\lambda\Phi_1L_\rho,\\
 c_\omega&=2k_P\lambda\Phi_1\eps/a.
 \end{aligned}
\end{equation}
If
\begin{equation}\label{eq:reg_trace_smallgain}
 a\sqrt\omega\ge4k_P\lambda\Phi_1K_0,
\end{equation}
then
\begin{equation}\label{eq:Lz}
 \begin{aligned}
 &\sup_{t\in I_T}\int_{t_0}^t
        \ee^{-\alpha(t-s)}|\dot z(s)|\,\dd s\le L_z,\\
 &L_z=2\left[H_{\rm in}+K_0U_\star
                   +\frac{L_H+K_0C_0}{\alpha}\right].
 \end{aligned}
\end{equation}
The bound is uniform in $T$, initial time, dither phase, and
frequency in the stated range.
\end{lemma}

\begin{proof}
The physical trace is known only to be continuous.
The existing mild solution already determines $H$, hence $Y_H$
and $z_{\rm f}$. We keep these functions fixed and solve the
resulting scalar Volterra equation in a class of absolutely
continuous trial traces. Uniqueness in the larger class of
continuous traces will then identify this constructed solution
with the physical trace. This establishes the regularity required
for differentiation with respect to time and the chain rule.
Lemma~\ref{lem:spatial} gives $z_{\rm f}'\in L^1_{\rm loc}$
and $\norm{z_{\rm f}'}_{L^1(t_0,t_0+T_1)}\to0$ as $T_1\downarrow0$.
Fix $R_d>0$ smaller than one quarter of the distance of
$z(t_0)+\rho_\omega(t_0)$ from the endpoints of the enlarged
trace interval. This distance is positive by
\eqref{eq:reg_stopped_region}. On a sufficiently short
$I=[t_0,t_0+T_1]$, the complete derivative ball
\[
 \begin{aligned}
 \mathcal D=\{\zeta\in W^{1,1}(I):\ &\zeta(t_0)=z(t_0),\norm{\dot\zeta}_{L^1(I)}\le R_d\}
 \end{aligned}
\]
has $\zeta+\rho_\omega$ in that interval. Its metric is the
$L^1$ distance of the derivatives. The common initial value
identifies this set with a closed ball in $L^1(I)$.
Scalar norms below are on $I$
unless another interval is indicated.

We implement the controller with $J_\zeta=\varphi(\zeta+\rho_\omega)$
from the actual initial controller states. Let
$f_\zeta=m\varphi'(\zeta+\rho_\omega)(\dot\zeta+\dot\rho_\omega)$.
The trial responses are absolutely continuous. Thus, \eqref{eq:filteridentity}
and the bias equation give
\begin{equation}\label{eq:reg_trial_controller}
 \begin{aligned}
 \dot{\widehat g}_\zeta&=-\lambda\widehat g_\zeta+f_\zeta,\\
 \dot U_\zeta&=(k_P\lambda-k_I)\widehat g_\zeta-k_Pf_\zeta.
 \end{aligned}
\end{equation}
Only the trial traces have been differentiated here. Absolute
continuity of the physical trace will follow by uniqueness.
For $d_\zeta=\norm{\dot\zeta_1-\dot\zeta_2}_1$,
the common initial value implies
$\norm{\zeta_1-\zeta_2}_\infty\le d_\zeta$. Therefore
\[
 \begin{aligned}
 \norm{f_{\zeta_1}-f_{\zeta_2}}_1&\le A_{T_1}d_\zeta,\\
 A_{T_1}&=\norm m_\infty
       [\Phi_1+\Phi_2(R_d+\norm{\dot\rho_\omega}_1)].
 \end{aligned}
\]
Variation of constants in \eqref{eq:reg_trial_controller} gives
$\norm{\widehat g_{\zeta_1}-\widehat g_{\zeta_2}}_\infty
\le A_{T_1}d_\zeta$ and hence
\begin{equation}\label{eq:reg_controller_difference}
 \begin{aligned}
 \norm{\dot U_{\zeta_1}-\dot U_{\zeta_2}}_1&\le L_{T_1}d_\zeta,\\
 L_{T_1}&=(|k_P\lambda-k_I|T_1+k_P)A_{T_1}.
 \end{aligned}
\end{equation}
This constant stays bounded as $T_1\downarrow0$ for fixed
$a,\omega,R_d$. Bounded $\varphi''$ suffices.

Next, we define $\mathcal T\zeta=z_{\rm f}+k*U_\zeta$. For absolutely
continuous $U$, inserting
$U(s)=U(t_0)+\int_{t_0}^s\dot U(r)\,\dd r$ and using Fubini gives
\begin{equation}\label{eq:reg_convolution_derivative}
 (k*U)'=k(t-t_0)U(t_0)+k*\dot U
 \quad\hbox{a.e.}
\end{equation}
All trial responses share $U(t_0)$. Thus
\eqref{eq:reg_controller_difference} gives
\begin{equation}\label{eq:reg_trace_contraction}
 \norm{(\mathcal T\zeta_1-\mathcal T\zeta_2)'}_1
 \le\norm k_{L^1(0,T_1)}L_{T_1}d_\zeta.
\end{equation}
The shared initial trace is gives the same initial
command. Therefore, $k(t-t_0)U(t_0)$ cancels when two responses are subtracted.
The remaining operator has norm proportional to
$\norm k_{L^1(0,T_1)}$, which tends to zero despite the kernel's
singularity at zero. Thus \eqref{eq:reg_trace_contraction} becomes
contractive without differentiation of $k$.
Also, \eqref{eq:reg_trial_controller} yields
\[
 \norm{\widehat g_\zeta}_\infty
 \le |\widehat g(t_0)|
    +\norm m_\infty\Phi_1(R_d+\norm{\dot\rho_\omega}_1),
\]
so $\norm{\dot U_\zeta}_1$ stays uniformly bounded over
$\mathcal D$. Equation~\eqref{eq:reg_convolution_derivative} implies
\[
 \norm{(\mathcal T\zeta)'}_1
 \le\norm{z_{\rm f}'}_1
       +\norm k_{L^1(0,T_1)}(|U(t_0)|+\norm{\dot U_\zeta}_1).
\]
This tends uniformly to zero. Since $\mathcal T$ preserves the
initial value, it maps $\mathcal D$ into itself and is a contraction
for small $T_1$. Its fixed point is absolutely continuous.

For continuous trial traces, variation of constants in the original
$p$ equation and integration of the bias equation give
\[
 \norm{U_{\zeta_1}-U_{\zeta_2}}_\infty
 \le\Lambda_{T_1}\norm{\zeta_1-\zeta_2}_\infty,
\]
where
\[
 \begin{aligned}
 \Lambda_{T_1}&=(k_P+k_IT_1)\Phi_1[\norm m_\infty
        +T_1(\lambda\norm m_\infty+\norm{\dot m}_\infty)].
 \end{aligned}
\]
The first factor accounts for the bias integration and proportional
term. The bracket bounds the change in $p+mJ$ using
$|J_{\zeta_1}-J_{\zeta_2}|\le\Phi_1|\zeta_1-\zeta_2|$.
For fixed $a,\omega$, $\norm k_{L^1(0,T_1)}\Lambda_{T_1}<1$
when $T_1$ is small. Hence the Volterra equation has at most one
continuous solution, identifying the fixed point with the physical trace.

At any restart time $t_1$, the bound
$\norm{e(t_1)}_{C^\gamma}\le E_\gamma$ gives
\[
 \begin{aligned}
 \int_{t_1}^{t_1+\tau_1}|z_{{\rm f},t_1}'(s)|\,\dd s
 &\le E_\gamma\left[\mu_0\tau_1
       +\frac{2c_A\kappa^{\gamma/2}}\gamma\tau_1^{\gamma/2}\right]+L_H\tau_1.
 \end{aligned}
\]
This tends to zero uniformly in $t_1$. On a compact interval
satisfying \eqref{eq:reg_stopped_region}, the trace has a positive
margin within its regularity interval and the controller coefficients
are bounded. The radius and contraction interval can therefore be
chosen uniformly over restart times. Finitely many restarts establish
absolute continuity on $I_T$.

We may now differentiate \eqref{eq:tracerepresentation}:
\begin{equation}\label{eq:reg_trace_derivative}
 \begin{aligned}
 \dot z&=\Cb AS(t-t_0)e(t_0)+\Cb Y_H'+k(t-t_0)U(t_0)+k*\dot U.
 \end{aligned}
\end{equation}
Using \eqref{eq:reg_trial_controller} for the physical trace,
$|\widehat g|\le G$, and the harmonic bounds gives
\begin{equation}\label{eq:reg_command_derivative_bound}
 |\dot U|\le C_0+c_\omega|\dot z|,
\end{equation}
with \eqref{eq:reg_trace_coefficients}. Denote the supremum in
\eqref{eq:Lz} by $V_T$. It is finite by absolute continuity.
For $e_\alpha(t)=\ee^{-\alpha t}$, commutation of convolutions gives
\[
 e_\alpha*(|k|*|\dot z|)
 =|k|*(e_\alpha*|\dot z|)\le K_0V_T.
\]
Equations~\eqref{eq:reg_time_bounds}, \eqref{eq:reg_trace_derivative},
and \eqref{eq:reg_command_derivative_bound} consequently yield
\[
 V_T\le H_{\rm in}+K_0U_\star
          +\frac{L_H+K_0C_0}{\alpha}+c_\omega K_0V_T.
\]
Condition~\eqref{eq:reg_trace_smallgain} gives $c_\omega K_0\le1/2$,
proving \eqref{eq:Lz}. The unweighted norm $K_0$ suffices even
when $\alpha\ge\mu_0$, and all constants are uniform as stated.
\end{proof}

%% file: sections/trajectory_comparison.tex
\section{Trajectory Comparison}\label{sec:comparison}
The regularity estimates allow us to compute the difference between
the transformed implementation and the averaged equations exactly.
We then apply the disturbance estimate of Lemma~\ref{lem:averaged}.
{The residuals collect everything the averaged target ignores which include the reaction evaluated on the lifted field, the second-harmonic and Taylor
remainders of the cost, and the finite-frequency corrections to the
demodulation. Each is proportional to a power of $\eps=\omega^{-1/2}$.}

\subsection{Exact residual equations}
On a nominal existence interval with $Z_\omega\in B_R(\Zs)$,
impose the hypotheses of Lemmas~\ref{lem:spatial} and~\ref{lem:trace}.
Then $z=\Cb x$ is absolutely continuous, so the scalar derivatives
below exist almost everywhere and their integral identities hold
throughout the interval.

The lifting removes the dither from the boundary condition while
leaving the reaction evaluated at $x+r_\omega$. The remaining
boundary command is $v-k_P\widehat g$. Using
$\widehat g=\xi+\Delta_\omega$ therefore writes the transformed
field and bias as the averaged target with disturbances
\begin{equation}\label{eq:disturbances}
\begin{aligned}
 d_w&=F(x+r_\omega)-F(x),\\
 d_u&=-k_P\Delta_\omega,\qquad d_v=-k_I\Delta_\omega.
\end{aligned}
\end{equation}
Here $d_u$ enters the boundary command. To compute the filter
residual, we define
\begin{equation}\label{eq:costresidual}
\begin{aligned}
 E_\varphi={}&\varphi(z+\rho_\omega)-\varphi(z)-a\eps g(x)\sin(\theta-\pi/4).
\end{aligned}
\end{equation}
By Taylor's formula and \eqref{eq:traceharmonic}, we obtain 
\begin{equation}\label{eq:costresidualbound}
 |E_\varphi|\le\Phi_1L_ea\eps^3
       +\tfrac12\Phi_2L_\rho^2a^2\eps^2.
\end{equation}
The first term comes from the error in the linear harmonic response. The second is the quadratic Taylor remainder of the cost.
Differentiation of \eqref{eq:coordinates} almost everywhere yields
\[
\begin{aligned}
 \dot\xi={}&-\lambda\xi-\lambda\Delta_\omega
 -\dot m[\varphi(z+\rho_\omega)-\varphi(z)]+m\varphi'(z)\dot z-\dot d_{2,\omega}.
\end{aligned}
\]
Inserting \eqref{eq:costresidual}, use the exact identities
\[
\begin{aligned}
 -\dot m\,a\eps g(x)\sin(\theta-\pi/4)
 &=\lambda g(x)[1-\sin(2\theta)],
 \end{aligned}
 \]
\[
\begin{aligned}
\dot d_{2,\omega}
 =\tfrac12\lambda\chi\eps^2\varphi''(z)\dot z\cos(2\theta)-\lambda g(x)\sin(2\theta).
\end{aligned}
\]
The first identity follows from
$2\cos(\theta+\pi/4)\sin(\theta-\pi/4)=\sin(2\theta)-1$.
The doubled-frequency terms cancel, leaving
\begin{equation}\label{eq:filterresidual}
\begin{aligned}
 \dot\xi&=-\lambda\xi+\lambda g(x)+d_\xi
\end{aligned}
\end{equation}
where
\[
 d_\xi=-\lambda\Delta_\omega-\dot mE_\varphi+\bigl[m\varphi'(z)
 -\tfrac12\lambda\chi\eps^2\varphi''(z)\cos(2\theta)\bigr]\dot z.
\]
Thus \eqref{eq:Deltabound} and \eqref{eq:costresidualbound} imply
\begin{equation}\label{eq:scalarresidual}
 |d_\xi|\le K_\xi[a\eps+\eps^2+(\eps/a+\eps^2)|\dot z|],
\end{equation}
where
\begin{equation}\label{eq:Kxi}
\begin{aligned}
 \hspace{-0.15in}K_\xi=\max\bigl\{&\lambda\Phi_2L_\rho^2,
 \lambda K_\Delta+2\lambda\Phi_1L_e,2\lambda\Phi_1,\ \tfrac12\lambda|\chi|\Phi_2\bigr\}.
\end{aligned}
\end{equation}
More explicitly, the four coefficients in \eqref{eq:Kxi} arise from
\[
\begin{aligned}
 |d_\xi|\le{}\lambda\Phi_2L_\rho^2a\eps &
       +(\lambda K_\Delta+2\lambda\Phi_1L_e)\eps^2\\
 &+\left(2\lambda\Phi_1\frac{\eps}{a}
        +\tfrac12\lambda|\chi|\Phi_2\eps^2\right)|\dot z|.
\end{aligned}
\]
The factor $\dot m$, whose magnitude grows with frequency, turns
the quadratic cost remainder into the term proportional to $a\eps$.
The coefficient of $|\dot z|$ decreases with frequency at fixed
$a$. Its effect is controlled by the weighted integral in
\eqref{eq:Lz}; no bounded time derivative at startup is required.
Only $C^2$ cost regularity is used. The other residuals satisfy
\begin{equation}\label{eq:otherresidualbounds}
\begin{aligned}
 \norm{d_w}_\infty&\le L_FL_ra\eps,\\
 |d_u|&\le k_PK_\Delta\eps^2,\\
 |d_v|&\le k_IK_\Delta\eps^2.
\end{aligned}
\end{equation}
The full $L_F$ is required because $d_w$ includes $V_0r_\omega$,
even when this linear reaction contribution supplies strong damping
to the averaged system. The smaller derivative variation $\ell(R)$
cannot replace it.

\subsection{A bound independent of the comparison interval}
In this section, we show that the comparison error bound is given by
\begin{equation}\label{eq:comparisonradius}
\begin{aligned}
 \Romega={}&G_B\frac\sigma\alpha k_PK_\Delta\eps^2
       +\frac M\alpha(L_FL_r+K_\xi)a\eps\\
 &+\frac M\alpha(k_IK_\Delta+K_\xi)\eps^2
       +MK_\xi(\eps/a+\eps^2)L_z.
\end{aligned}
\end{equation}
We note that every constant on the right is fixed independently of the comparison interval. For fixed $a>0$, $\Romega=O(\omega^{-1/2})$.
\begin{lemma}\label{lem:comparison}
Under Assumptions~\ref{ass:regularity}--\ref{ass:margin} and
\eqref{eq:gaincondition}, fix $0<a\le a_0$ and impose
\eqref{eq:initialfamily}. Suppose $\omega\ge\omega_{\rm lift}$,
$L_ra\eps\le\eta_\star$, and \eqref{eq:reg_trace_smallgain} hold.
For transformed nominal and averaged solutions remaining in
$B_R(\Zs)$ throughout $[t_0,t]$,
\begin{equation}\label{eq:comparison_general}
\begin{aligned}
 \norm{Z_\omega(t)-\bar Z(t)}_\E
 \le{}&M\ee^{-\alpha(t-t_0)}
       \norm{Z_\omega(t_0)-\bar Z(t_0)}_\E\\
 &+\Romega
\end{aligned}
\end{equation}
With matched initial states, the following bounds hold
\begin{equation}\label{eq:stoppedcomparison}
\begin{aligned}
 \norm{Z_\omega(t)-\bar Z(t)}_\E&\le\Romega,\\
 \norm{X_\omega(t)-\bar Z(t)}_\E&\le\Romega+\domega.
\end{aligned}
\end{equation}
The bounds are independent of the interval length.
\end{lemma}
\begin{proof}
The exact residual equations \eqref{eq:disturbances} and
\eqref{eq:filterresidual} meet the hypotheses of
Lemma~\ref{lem:averaged}. The field and boundary residuals are
continuous, and the filter residual is locally integrable by
Lemma~\ref{lem:trace}. The boundary contribution in
\eqref{eq:master} is at most $G_B\sigma k_PK_\Delta\eps^2/\alpha$.
For the internal terms, \eqref{eq:scalarresidual} and
\eqref{eq:otherresidualbounds} give
\[
\begin{aligned}
 &M\int_{t_0}^t\ee^{-\alpha(t-s)}
   [\norm{d_w(s)}_\infty+|d_v(s)|+|d_\xi(s)|]\,\dd s\\
 &\quad\le\frac M\alpha\bigl[(L_FL_r+K_\xi)a\eps
                         +(k_IK_\Delta+K_\xi)\eps^2\bigr]\\
 &\qquad+MK_\xi(\eps/a+\eps^2)
       \int_{t_0}^t\ee^{-\alpha(t-s)}|\dot z(s)|\,\dd s.
\end{aligned}
\]
The last integral is at most $L_z$ by \eqref{eq:Lz}.
Adding the boundary contribution yields \eqref{eq:comparisonradius}
and proves \eqref{eq:comparison_general}.
Matched initialization removes the initial-error term, and
\eqref{eq:coordbound} gives \eqref{eq:stoppedcomparison}.
\end{proof}

This is a conditional comparison on the region where the estimates
apply. The next section combines it with averaged decay to show that
small initial errors and sufficiently high frequencies keep the
implemented solution in that region for all time.

%% file: sections/practical_stability.tex
\section{Practical Stability}\label{sec:practical}
\begin{table}[t]
\caption{{Main constants of the analysis and where they are defined.}}
\label{tab:constants}
\centering
{%
\footnotesize
\setlength{\tabcolsep}{3pt}
\renewcommand{\arraystretch}{1.12}
\begin{tabular}{@{}lll@{}}
\toprule
Symbol & Meaning & Defined in\\
\midrule
$\eps$, $\omega_{\rm lift}$ & $\eps=\omega^{-1/2}$; minimal lifting frequency
 & Secs.~\ref{sec:controller}, \ref{sec:harmonic}\\
$M$, $\sigma$ & overshoot and decay rate of the linear loop
 & \eqref{eq:referenceconstants}, \eqref{eq:M}\\
$G_B$ & gain from a boundary command to the state & \eqref{eq:GBdef}\\
$\ell(R)$, $\alpha$ & nonlinear margin; $\alpha=\sigma-M\ell(R)$
 & \eqref{eq:ell}\\
$\Phi_1$, $\Phi_2$ & bounds on $\varphi'$, $\varphi''$ near the optimum
 & after \eqref{eq:ell}\\
$L_r$, $L_\rho$, $L_e$ & lifted field, its trace, harmonic remainder
 & \eqref{eq:liftconstants}\\
$K_\Delta$, $\domega$ & filter correction; coordinate error
 & \eqref{eq:Deltabound}, \eqref{eq:delta}\\
$K_0$ & $L^1$ norm of the trace kernel & \eqref{eq:reg_K0}\\
$L_z$ & weighted bound on the trace derivative & \eqref{eq:Lz}\\
$K_\xi$ & filter residual coefficient & \eqref{eq:Kxi}\\
$\Romega$ & implemented-vs-averaged comparison error
 & \eqref{eq:comparisonradius}\\
$\bomega$, $r_\star$ & ultimate bound; admissible initial radius
 & \eqref{eq:practicalconstants}\\
\bottomrule
\end{tabular}}
\end{table}
{Table~\ref{tab:constants} collects the constants that enter the
final estimates.}
By the developments above, the averaged decay and the comparison bound now exclude departure
from the region where the estimates were obtained. With
\eqref{eq:delta} and \eqref{eq:comparisonradius}, we let
\begin{equation}\label{eq:practicalconstants}
 \begin{aligned}
 \bomega&=\Romega+(M+1)\domega,\\
  r_\star&=R/(8M).
  \end{aligned}
\end{equation}
The term $\domega$ bounds the coordinate change, while $\Romega$
bounds the error dynamics. The additional $M\domega$ in
$\bomega$ accounts for the change in the initial condition when
applying the averaged decay estimate. This initial contribution
decays, which explains the smaller ultimate radius below.
We impose the sufficient frequency conditions
\begin{equation}\label{eq:frequency}
\begin{aligned}
 \omega&\ge\omega_{\rm lift},\qquad
 a\sqrt\omega\ge4k_P\lambda\Phi_1K_0,\\
 L_ra/\sqrt\omega&\le\eta_\star,\qquad \bomega\le R/8.
\end{aligned}
\end{equation}
For each fixed $a>0$, they hold for all frequencies above a finite
threshold, since $\domega,\Romega\to0$.

\begin{theorem}\label{thm:unified}
Suppose Assumptions~\ref{ass:regularity}--\ref{ass:margin} and
\eqref{eq:gaincondition} hold. Fix $a\in(0,a_0]$ and the
initial-family bound $H_0$. Choose $\omega$ satisfying
\eqref{eq:frequency}. For every $t_0\ge0$, every dither phase,
and every initial state satisfying \eqref{eq:initialfamily} and
$E_\omega(t_0)\le r_\star$, the implemented loop has a unique
global mild solution. Initialize the averaged comparison by
\begin{equation}\label{eq:initialization}
 \bar Z(t_0)=Z_\omega(t_0)
\end{equation}
using \eqref{eq:coordinates}. Then, for all $t\ge t_0$,
\begin{equation}\label{eq:joint}
\begin{aligned}
 \norm{Z_\omega(t)-\bar Z(t)}_\E&\le\Romega,\\
 \norm{X_\omega(t)-\bar Z(t)}_\E&\le\Romega+\domega,\\
 E_\omega(t)&\le M\ee^{-\alpha(t-t_0)}E_\omega(t_0)+\bomega.
\end{aligned}
\end{equation}
Moreover,
\begin{equation}\label{eq:ultimate}
\begin{aligned}
 \limsup_{t\to\infty}E_\omega(t)&\le\Romega+\domega,\\
 \limsup_{t\to\infty}|J(t)-\varphi(z^\star)|
 &\le\tfrac12\Phi_2(\Romega+L_ra\eps)^2.
\end{aligned}
\end{equation}
The state and cost bounds have orders $\omega^{-1/2}$ and
$\omega^{-1}$ at fixed $a$. They are uniform in initial time,
phase, admitted initial state, and frequency above the selected
threshold.
\end{theorem}
\begin{proof}
For fixed $a,\omega$, the original mild field equation and the
controller integral equations are locally contractive: $\KB$ is
locally integrable and all state nonlinearities are locally
Lipschitz. This gives local existence without differentiating $J$.

First justify the initial coordinate bound. Condition
\eqref{eq:frequency} implies
$\domega\le R/[8(M+1)]$ and $\Romega\le R/8$, so
\[
 \norm{x(t_0)-w^\star}_\infty
 \le r_\star+L_ra\eps\le r_\star+\domega<R.
\]
The lift stays within the enlarged regularity region. Applying
Lemma~\ref{lem:coordinates} therefore gives
\begin{equation}\label{eq:initialmargin}
 M\norm{Z_\omega(t_0)-\Zs}_\E
       \le M(r_\star+\domega)<R/4.
\end{equation}
By \eqref{eq:initialmargin}, Corollary~\ref{cor:target} applies to
the matched initial state \eqref{eq:initialization}. Its decay estimate
\eqref{eq:targetbound} keeps $\bar Z$ in $B_{R/4}(\Zs)$.
Before any nominal first exit from $B_R(\Zs)$,
\eqref{eq:stoppedcomparison} gives
\[
 \norm{Z_\omega(t)-\Zs}_\E<R/4+\Romega\le3R/8.
\]
By continuity, any such exit is excluded.

The inverse transformation bounds $w,v,\widehat g$ and
$p=\widehat g-m\varphi(\Cb w)$ on every finite interval.
Since $3R/8+\eta_\star<R+\eta_\star$, the physical field remains a positive distance inside the regularity region.
At a finite maximal time, the bounded internal perturbation and the strong
semigroup continuity yield an endpoint limit of the distributed
convolution. The boundary convolution has an endpoint limit by
local integrability of $\KB$ and boundedness of its command.
The controller states also have limits. These facts allow one to restart the
local mild problem while excluding finite maximal existence.
Boundedness of the cost is used here only for continuation,
not in the error constants.

The global comparison, averaged decay, and coordinate estimate
now imply
\begin{equation}\label{eq:link}
\begin{aligned}
 E_\omega(t)
 &\le M\ee^{-\alpha(t-t_0)}[E_\omega(t_0)+\domega]
            +\Romega+\domega\\
 &\le M\ee^{-\alpha(t-t_0)}E_\omega(t_0)+\bomega.
\end{aligned}
\end{equation}
Estimate~\eqref{eq:link} proves \eqref{eq:joint} and the state part of
\eqref{eq:ultimate}. For the field alone,
$\limsup\norm{w-w^\star}_\infty\le\Romega+L_ra\eps$.
Since $\varphi'(z^\star)=0$ and $\norm{\Cb}=1$, Taylor's formula
bounds the cost gap by $\Phi_2\norm{w-w^\star}_\infty^2/2$,
proving the remaining assertion. All estimates use elapsed time
and phase-uniform bounds. Their constants are fixed over the
prescribed initial family for each chosen $a>0$.
\end{proof}

\begin{corollary}\label{cor:practical}
Under Theorem~\ref{thm:unified}, fix $\eta>0$ and strengthen
\eqref{eq:frequency} to $\bomega\le\min\{R/8,\eta/2\}$.
Every admitted state with initial error at most
$\min\{r_\star,\eta/(2M)\}$ satisfies $E_\omega(t)\le\eta$
for all $t\ge t_0$. Every admitted state with error at most
$r_\star$ satisfies that accuracy after
\begin{equation}\label{eq:time}
 t-t_0\ge T_\eta:=\alpha^{-1}
       \max\{0,\log(2Mr_\star/\eta)\}.
\end{equation}
\end{corollary}
\begin{proof}
For the smaller initial radius, or after the delay in \eqref{eq:time},
the decaying term in \eqref{eq:joint} is at most $\eta/2$. Adding
$\bomega\le\eta/2$ proves both assertions.
\end{proof}

The matched averaged initial state depends on frequency and
phase. The result gives practical regulation with the dither
active. Invariance of the ultimate ball and existence or
uniqueness of a periodic orbit are not asserted.
{In the thermal setting, Theorem~\ref{thm:unified} says that a
plate started near its optimal temperature profile, with the filter
initialized by \eqref{eq:filterinitialization}, stays near that profile
for all time, that the element temperature settles to within a band of
width proportional to $\omega^{-1/2}$ around its optimum, and that the
residual loss of transmission is of order $\omega^{-1}$.}

%% file: sections/frequency_and_scope.tex
\section{Frequency Selection and Scope}\label{sec:frequency}

By \eqref{eq:practicalconstants}, we obtain
\begin{equation}\label{eq:frequencyerror}
 \bomega=\Gamma_1(a)/\sqrt\omega+\Gamma_2/\omega,
\end{equation}
where
\begin{equation}\label{eq:gammas}
\begin{aligned}
 \Gamma_1(a)
 &=a\left[\frac M\alpha(L_FL_r+K_\xi)\right.\\
&\hspace{0.75in}\left.+(M+1)L_r\right]+\frac{MK_\xi L_z}{a},\\
 \Gamma_2
 &=G_B\frac\sigma\alpha k_PK_\Delta
       +\frac M\alpha(k_IK_\Delta+K_\xi)\\
 &\quad+MK_\xi L_z+(M+1)K_\Delta.
\end{aligned}
\end{equation}
For $0<s\le R/8$, solving \eqref{eq:frequencyerror} and imposing
the remaining conditions of \eqref{eq:frequency} gives
\begin{equation}\label{eq:completefrequency}
\begin{aligned}
 \sqrt\omega\ge\max\biggl\{&\sqrt{\omega_{\rm lift}},\,
 \frac{4k_P\lambda\Phi_1K_0}{a},\,\frac{L_ra}{\eta_\star},\\
 &\frac{\Gamma_1(a)+\sqrt{\Gamma_1(a)^2+4\Gamma_2s}}{2s}
 \biggr\}.
\end{aligned}
\end{equation}
Then $\bomega\le s$. We use $s=R/8$ for Theorem~\ref{thm:unified}
or $s=\min\{R/8,\eta/2\}$ for Corollary~\ref{cor:practical}.
The constants in \eqref{eq:gammas} use  fixed $a_0$, leaving 
terms proportional to $a$ and $1/a$. {The terms proportional to $a$
come from the lifted field and from the reaction evaluated on it. The
term proportional to $1/a$ comes from the demodulation gain, which
grows as the dither amplitude shrinks.} Smaller dither amplitudes do
not necessarily improve the bound. The controller does not evaluate
\eqref{eq:completefrequency}, but certification requires plant and
cost-derivative bounds.

For $F(w)=-Kw-\beta w^3$, $d=0$, and $\varphi(z)=z^2/2$, with
$K,\beta>0$, the optimum is $w^\star=u^\star=0$. Choosing
$V_0=-K$ gives $\mu_{\rm ref}=\mu_0+K$ and $\ell(R)=3\beta R^2$.
Thus $R\le\sqrt{\sigma/(6M\beta)}$ ensures $\alpha\ge\sigma/2$
for every fixed $K$. The full bound
$L_F=K+3\beta(R+\eta_\star)^2$ still enters the approximation
error, so absorbing damping need not reduce the sufficient frequency.
Transverse slope variation remains subject to
Assumption~\ref{ass:margin}.

{The universal $\pi/4$ phase is a property of the heated edge.
Moving the sensor into the plate adds a frequency-dependent delay and
attenuation, which is why the design assumes a sensor on the actuated
boundary.} For the linear diffusion plant, an interior midpoint sensor
at $(d_s,1/2)$, $0<d_s<1$, has
\[
\begin{aligned}
 \frac{H_{d_s}(i\omega)}{H_b(i\omega)}
 &=\frac{\cosh[\gamma_0(1-d_s)]}{\cosh\gamma_0}\\
 &=\exp\!\left[-(1+i)d_s\sqrt{\frac{\omega}{2\kappa}}\right]
       [1+o(1)].
\end{aligned}
\]
Hence the boundary phase in \eqref{eq:asymptotic}, after separating
the signed gain, cannot be used unchanged for interior sensing.
Both its extra attenuation and phase depend on frequency.

For a prescribed physical cost trace, bounded measurement noise
$n$ and zero initial perturbation of $p$ produce
\[
\begin{aligned}
 q_n(t)&=m(t)n(t)\\
 &\quad-\int_{t_0}^t\ee^{-\lambda(t-s)}
       [\lambda m(s)+\dot m(s)]n(s)\,\dd s,
\end{aligned}
\]
and therefore, almost everywhere,
\begin{equation}\label{eq:noisefilter}
 |q_n(t)|\le\left(\frac{4\lambda\eps}{a}
              +\frac2{a\eps}\right)\norm n_{\infty,[t_0,t]}.
\end{equation}
The growing term in \eqref{eq:noisefilter} concerns noise near the
dither frequency. For $n=N\sin(\theta+\psi)$, the steady filter
mean is $N\sqrt\omega\sin(\pi/4-\psi)/a$.
Constant offsets instead cancel under consistent initialization.
These are prescribed-signal filter calculations, not noisy
closed-loop guarantees. The first-order cost ripple also vanishes
at the optimum because $\varphi'(z^\star)=0$.

The theorem concerns local regulation to an attainable output-cost
minimum. Global optimization, input saturation, measurement delays,
and additional actuator or sensor dynamics are outside its scope.

%% file: sections/physical_example.tex
\section{Thermal Resonance Alignment}\label{sec:physicalexample}

{We now return to the plate of Fig.~\ref{fig:geometry}(a) and give
it numbers.} Thermally tunable filters motivate the local optical
model \cite{Domash2004,HohlfeldZappe2007}.

\begin{figure*}[!t]
\centering
\includegraphics[width=0.62\textwidth]{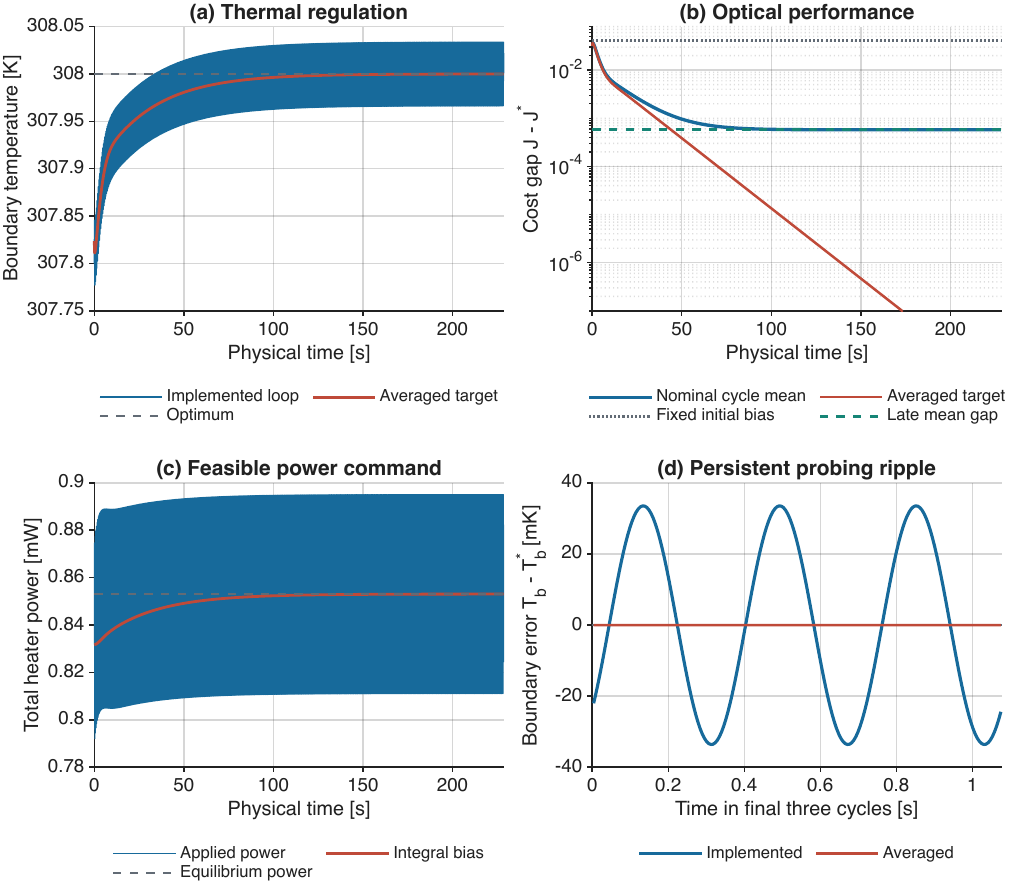}
\caption{Controller performance at $\omega=80$: (a) boundary temperature,
(b) optical cost gap, (c) applied power and integral bias, and
(d) {boundary-temperature error over the final three dither
cycles}. The nominal cost in
(b) is averaged over complete dither cycles. The baseline holds the
initial heater power, and the averaged target uses matched coordinates.}
\label{fig:caseloop}
\end{figure*}

\begin{figure}[!t]
\centering
\includegraphics[width=0.9\columnwidth,trim=30 84 17 78,clip]{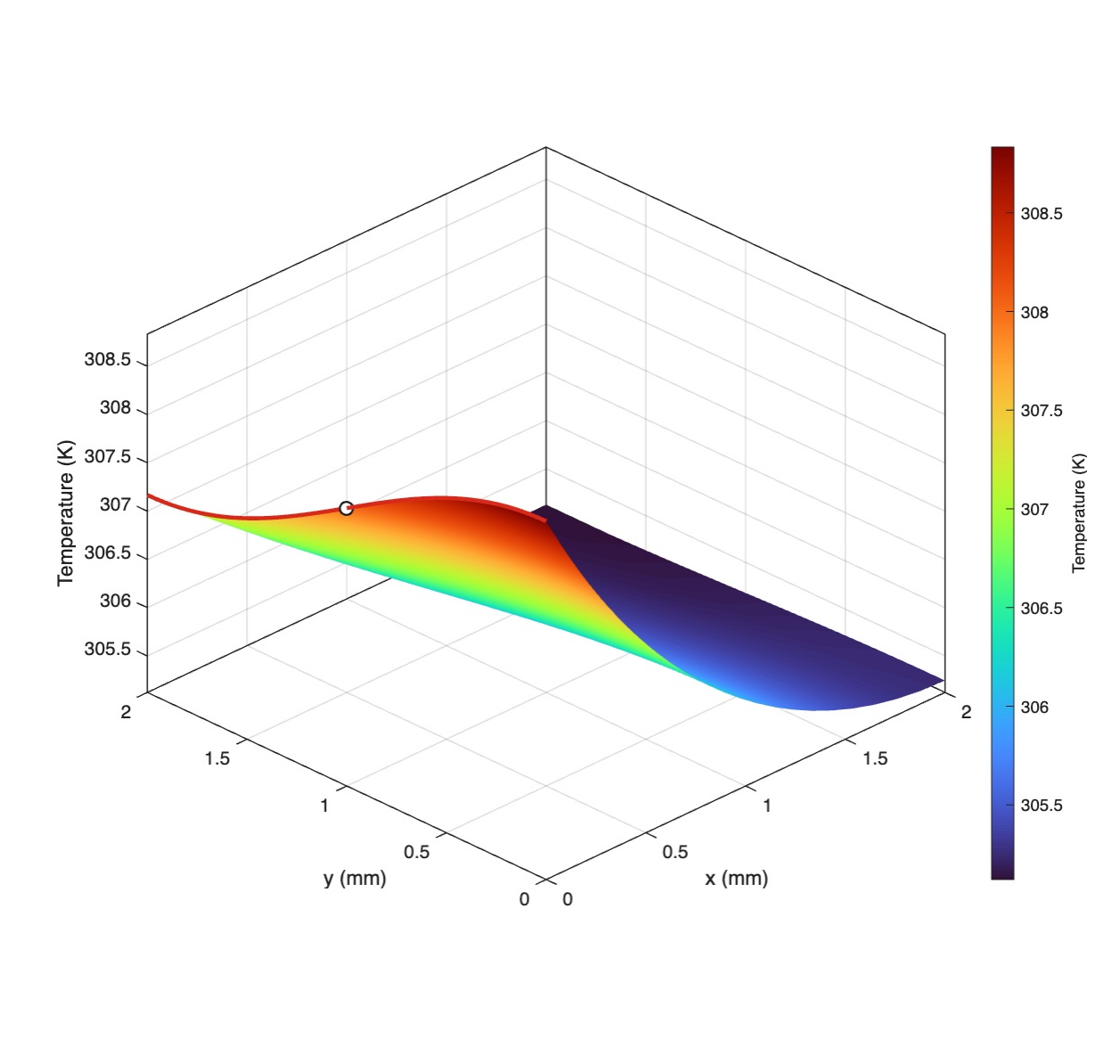}
\caption{{Stationary optimal temperature field $T^\star$ over the
plate. The red line marks the heated edge, along which the nonuniform
heater profile $q$ shapes the temperature, and the circle marks the
optical element, held at $308\,\mathrm K$.}}
\label{fig:casefields}
\end{figure}

\begin{figure*}[!t]
\centering
\includegraphics[width=\textwidth,trim=61 100 31 80,clip]{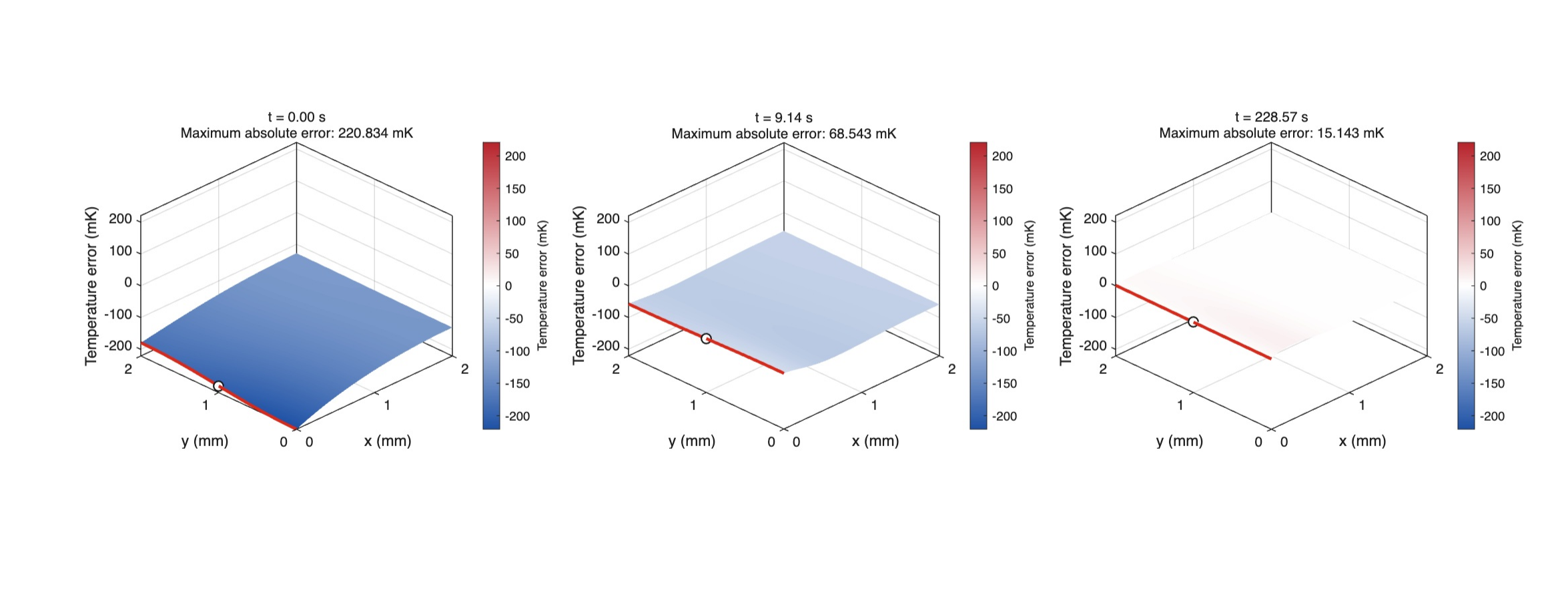}
\caption{{Error $T-T^\star$ between the implemented and the optimal
temperature fields at physical times $0$, $9.14$, and
$228.57\,\mathrm s$ (left to right); the maximum absolute error falls
from $221$ to $68.5$ to $15.1\,\mathrm{mK}$. The red line marks the
heated edge and the circle the optical element. The snapshots display
the field component, with neither the bias nor the sensitivity
estimate.}}
\label{fig:caseerror}
\end{figure*}

\subsection{Thermal model and scaling}

Let $T(X,Y,t_{\rm phys})$ be the temperature of a plate occupying
$(X,Y)\in(0,L)^2$. The thin-plate approximation treats temperature
as uniform through its thickness $d_p$. In-plane conduction and
linear cooling through the two broad faces then give
\begin{equation}\label{eq:casephysical}
\begin{aligned}
 C_vd_p\frac{\partial T}{\partial t_{\rm phys}}
   &=k_{\rm th}d_p\Delta_{X,Y}T-2h_f(T-T_a),\\
 -k_{\rm th}d_p T_X(0,Y,t_{\rm phys})
   &=\frac{P(t_{\rm phys})}{L}q(Y/L),\\
 T_X(L,Y,t_{\rm phys})&=0,\\
 T_Y(X,0,t_{\rm phys})&=T_Y(X,L,t_{\rm phys})=0.
\end{aligned}
\end{equation}
Here $C_v$ is volumetric heat capacity, $k_{\rm th}$ is thermal
conductivity, and $h_f$ is the heat-transfer coefficient of each
broad face. The factor two accounts for cooling from both faces.
The ambient temperature is $T_a$, and the three unactuated
lateral edges are insulated.

The scalar command $P$ is distributed along the heated edge by
$q(s_2)=1+\vartheta\cos(\pi s_2)$. The normalization
$\int_0^1q(s_2)\,ds_2=1$ ensures that integrating the imposed
edge flux in \eqref{eq:casephysical} gives the total heater
power $P$. The nonuniform profile produces a spatially varying
temperature along the heated edge, although the optical element
measures performance at its midpoint.

Table~\ref{tab:casephysical} lists the illustrative physical and
optical parameters used in the reported computation. The optical
half-width provides the temperature scale $\Theta$.

\begin{table}[t]
\caption{Physical and optical parameters of the thermal example.}
\label{tab:casephysical}
\centering
\small
\setlength{\tabcolsep}{3pt}
\renewcommand{\arraystretch}{1.12}
\begin{tabular}{@{}lll@{}}
\hline
Quantity & Symbol & Value\\
\hline
Plate side length & $L$ & $2\,\mathrm{mm}$\\
Thickness & $d_p$ & $0.1\,\mathrm{mm}$\\
Conductivity & $k_{\rm th}$ & $1.4\,\mathrm{W\,m^{-1}K^{-1}}$\\
Volumetric heat capacity & $C_v$ &
 $1.6\times10^6\,\mathrm{J\,m^{-3}K^{-1}}$\\
Each-face heat transfer & $h_f$ &
 $17.5\,\mathrm{W\,m^{-2}K^{-1}}$\\
Ambient temperature & $T_a$ & $300\,\mathrm K$\\
Temperature scale & $\Theta$ & $1\,\mathrm K$\\
Heater nonuniformity & $\vartheta$ & $0.45$\\
Peak optical transmission & $A_{\rm opt}$ & $0.98$\\
Alignment temperature & $T_b^\star$ & $308\,\mathrm K$\\
\hline
\end{tabular}
\end{table}

Introduce $s=(X/L,Y/L)$ and normalize temperature, time, and
power as follows
\begin{equation}\label{eq:casescaling}
\begin{aligned}
 w&=(T-T_a)/\Theta, & t&=t_{\rm phys}/t_D,\\
 t_D&=C_vL^2/k_{\rm th}, & u&=P/P_s,\\
 P_s&=k_{\rm th}\Theta d_p.
\end{aligned}
\end{equation}
The diffusion time $t_D$ makes the normalized conduction
coefficient unity. The power scale $P_s$ similarly makes the
boundary input gain unity, so the edge condition becomes
$-w_{s_1}=q\,u$. Substituting \eqref{eq:casescaling} into
\eqref{eq:casephysical} gives \eqref{eq:plant} with
\[
 \kappa=b=1,\qquad
 \mu_0=\frac{2h_fL^2}{k_{\rm th}d_p}=1,\qquad F=d=0.
\]
Numerically, $t_D=4.5714\,\mathrm s$ and
$P_s=0.14\,\mathrm{mW}$. Thus a dimensionless input specifies
total heater power, and the physical midpoint temperature is
$T_b=T_a+\Theta\Cb w$. These relations connect the PDE variables
and the physical axes in the performance plots.

\subsection{Optical score and numerical procedure}

Let $P_{\rm in}$ be fixed incident optical power and $P_{\rm tr}$
the transmitted power. The local Lorentzian model is
\[
 \frac{P_{\rm tr}}{P_{\rm in}}
 =\frac{A_{\rm opt}}{1+(\Cb w-z^\star)^2}.
\]
The temperature half-width equals $\Theta$, so the detuning
$\Cb w-z^\star$ is dimensionless. The reciprocal-power score
{plotted in Fig.~\ref{fig:geometry}(b)} is
\begin{equation}\label{eq:casecost}
 J=\frac{P_{\rm in}}{P_{\rm tr}}-1
 =\frac1{A_{\rm opt}}-1+\frac{(\Cb w-z^\star)^2}{A_{\rm opt}},
\end{equation}
with $A_{\rm opt}=0.98$. The chosen optimum $z^\star=8$
corresponds to $T_b^\star=308\,\mathrm K$. Only $J$ enters
the feedback. Neither the alignment temperature nor
\eqref{eq:casecost} is evaluated by the controller.
The stationary midpoint gain is
$G_0=\coth(1)$, so $u^\star=z^\star/G_0$ is attainable.
The stationary field $w^\star$ is spatially nonuniform because
heat is applied at one edge and lost through both faces
{(Fig.~\ref{fig:casefields})}.
Regulation therefore concerns this attainable field and its
required input bias. Since
$\varphi''=2/A_{\rm opt}$, $V_0=0$ gives $\ell(R)=0$, and the
averaged error dynamics equal the linear system of
Proposition~\ref{prop:linear}.

Use \eqref{eq:piarchitecture}--\eqref{eq:nominal} with
$\lambda=2$, $k_P=0.35$, $k_I=0.10$, $a=0.30$, and
$\omega=80$ ($2.7852\,\mathrm{Hz}$). Set $t_0=\theta_0=0$.
The initial field is the discrete stationary solution with
$\Cb w_0=7.8$, $v_0$ its corresponding bias, and
$p_0=-m(0)J(0)$ as in \eqref{eq:filterinitialization}, giving
$\widehat g(0)=0$.
The target uses \eqref{eq:initialization}.
The horizon $0\le t\le50$ is $228.6\,\mathrm s$.

The Python/SciPy computation uses continuous P1 elements
(2,009 nodes, 3,840 triangles graded toward the heated edge) and
implicit midpoint integration with 64 steps per dither period.
Time and mesh refinement change the boundary temperature by at most
$0.052$ and $0.023\,\mathrm{mK}$, respectively, over $0\le t\le5$.

\subsection{Controller performance}

Figure~\ref{fig:caseloop} shows regulation from $307.8$ toward
$308\,\mathrm K$. The bias moves from $0.8317$ to
$0.8531\,\mathrm{mW}$, while the applied power remains between
$0.7920$ and $0.8951\,\mathrm{mW}$ without imposed saturation.
Over $40\le t\le50$, the peak boundary error is
$33.8\,\mathrm{mK}$ and the mean cost gap is $5.76\times10^{-4}$,
a $98.6\%$ reduction from the fixed-bias gap $4.08\times10^{-2}$.
The leading ripple $\Theta a\chi/\sqrt\omega=33.5\,\mathrm{mK}$
contributes $a^2\chi^2/(2A_{\rm opt}\omega)=5.74\times10^{-4}$
to the mean gap in \eqref{eq:casecost}, consistent with the observed
residual loss.

{Figures~\ref{fig:casefields} and~\ref{fig:caseerror} show what
the controller does to the plate as a whole. The optimal field
$T^\star$ in Fig.~\ref{fig:casefields} is hottest along the heated
edge, where the nonuniform heater profile $q$ leaves its imprint, and
decays across the plate as heat escapes through the faces; holding the
element at $308\,\mathrm K$ means maintaining this entire profile, and
the controller does so without ever measuring it.
Figure~\ref{fig:caseerror} follows the error $T-T^\star$.} The
initial peak field error is $221\,\mathrm{mK}$ near the heated edge.
{It falls to $68.5\,\mathrm{mK}$ after two diffusion times
($9.14\,\mathrm s$) and to $15.1\,\mathrm{mK}$ at the end of the run
($228.57\,\mathrm s$), inside the $\pm34\,\mathrm{mK}$ ripple band of
Fig.~\ref{fig:caseloop}(d); what remains is the thin boundary layer
of the persistent dither along the heated edge.} The snapshots
illustrate the field component at selected times, not the full error
$E_\omega$ or an all-time numerical bound. 

Increasing $\omega$ from $20$ to $320$ reduces the late peak
boundary error from $67.7$ to $17.0\,\mathrm{mK}$ and mean cost gap
from $2.33\times10^{-3}$ to $1.44\times10^{-4}$. The slopes
$-0.50$ and $-1.00$ agree with the predicted frequency dependence
over this range. The plant, gains, amplitude, physical initial field,
bias, and rule $\widehat g(0)=0$ are fixed. The stored $p_0$ and
matched averaged initial state vary with frequency.

Condition \eqref{eq:strictification} permits $\delta_L<0.78125$.
With $\delta_L=0.78$, Proposition~\ref{prop:linear} gives
$\sigma=0.026$ and $M\approx134$, compared with the averaged
rightmost eigenvalue $-0.153$. However, $E_\omega(0)=0.373$
requires $R\ge8ME_\omega(0)\approx400$. Since
$\Phi_1\ge2R/A_{\rm opt}$ and $K_0=\coth(1)$, the condition
$a\sqrt\omega\ge4k_P\lambda\Phi_1K_0$ in \eqref{eq:frequency}
alone then requires $\omega>10^8$. These runs lie outside the
sufficient certificate of Theorem~\ref{thm:unified} and provide
performance evidence rather than verification of its hypotheses.
{The gap between the certificate and the observed behavior
reflects the conservatism of the constants $M$ and $\Phi_1$. The
frequency trends predicted by the theorem are the ones observed.}

%% file: sections/conclusion.tex
\section{Conclusion}

The proposed dual-mode controller uses only boundary cost measurements
to regulate a diffusion process near an attainable optimal equilibrium
without online PDE compensation {. It does so by relying on a
property of diffusion rather than on a model of the plant: the $\pi/4$
phase lag of a periodically heated surface fixes the demodulation, and
the unknown signed gain disappears into a positive loop gain}.
Stability of the averaged system and trajectory estimates valid for all
time establish local practical exponential regulation, with residual
state and cost errors decreasing as the dither frequency increases.
The thermal example illustrates the controller's performance, with
temperature and cost errors consistent with the predicted frequency
dependence {, and a temperature field that settles on the optimal
profile across the whole plate although only the element's transmission
is measured}.
Future work will address less conservative frequency conditions,
measurement noise, interior sensing, and additional device dynamics.